# Magic Tours of Knight in 4x4x4 Cube

Awani Kumar. 2/11, Vijayant Khand, Gomti Nagar, Lucknow 226010. INDIA

Mail id: awanieva@gmail.com

**Abstract**

*There are 216 magic tours of knight in 4x4x4 cube.*

**Introduction**: Knight is a unique piece. Its crooked move has fascinated not only chess players but also mathematics enthusiasts and some great mathematical minds, such as; Raimond de Montmort, Abraham de Moivre, and Leonhard Euler have delved into the knight's tour puzzle. According to Jelliss [1], oldest record of knight's tour on 8x8 board dates back to around 840 AD. For centuries, the traditional study of knight's tour was mostly confined to square boards. Vandermonde [2], a mathematician, was the first in constructing a 3-dimensional knight's tour in a $4 \times 4 \times 4$ cube in the year 1771. Later, Schubert [3], Gibbins [4], Stewart [5], Jelliss [6] and Petkovic [7] have delved into the knight's tour puzzle in 3-dimensional space. How many magic tours are there in 4x4x4 cube? The author plans to look at this question.

**Details**: Magic tour of knight in 4x4x4 cube is a pretty recent field of study. Stertenbrink [8] was the first person to discover a magic tour of knight in 4 x 4 x 4 cube in the year 2003. The sum of its space diagonals (Fig.11) is also equal to the magic constant 130. It is a 4-fold cyclic tour consisting of only squares and diamonds pattern. Awani Kumar [9] [10] [11] has looked into the possibilities of knight's tour in cubes and cuboids, having magic properties. There are 68 primary magic tours in 4x4x4 cube out of which, 54 are closed (or reentrant) and 14 are open magic tours. Considering their cyclic and reverse properties, there are 216 magic tours (188 closed tours and 28 open tours) as per Frenicle's method of classification. 48 tours are diagonally magic, that is, all the four space diagonals sum up to 130. All the tours are either squares and diamonds type or irregular tours. William Beverley constructed the first magic tour of knight on 8x8 board in the year 1848 using Beverley pattern. Surprisingly, here none of the tours have Beverley pattern. We can consider 4x4x4 cube as composed of 8 smaller 2x2x2 cubes and here, sum of each smaller cube is 260, that is, $1/8^{th}$ of the 4x4x4 cube.

**Conclusion**: The author has constructed magic knight tours in 8x8x8 and 12x12x12 cubes but the one in 6x6x6 (or 10x10x10 and such singly-even) cube is still to be discovered. Readers are requested to look for them. Awani Kumar [12] has extended knight's tour in higher dimensions and has enumerated over 200 magic tours in 4x4x4x4 hypercube. Readers are also requested to look for new ones there.

**Acknowledgement**: The author is grateful to Atulit Kumar for his help in preparation of this article.

# Magic Tours of Knight in 4x4x4 Cube

| 1st Layer | | | | | 2nd Layer | | | | | 3rd Layer | | | | | 4th Layer | | | |
|---|---|---|---|---|---|---|---|---|---|---|---|---|---|---|---|---|---|---|
| 1 | 20 | 47 | 62 | | 48 | 61 | 2 | 19 | | 29 | 16 | 51 | 34 | | 52 | 33 | 30 | 15 |
| 28 | 9 | 54 | 39 | | 53 | 40 | 27 | 10 | | 8 | 21 | 42 | 59 | | 41 | 60 | 7 | 22 |
| 45 | 64 | 3 | 18 | | 4 | 17 | 46 | 63 | | 49 | 36 | 31 | 14 | | 32 | 13 | 50 | 35 |
| 56 | 37 | 26 | 11 | | 25 | 12 | 55 | 38 | | 44 | 57 | 6 | 23 | | 5 | 24 | 43 | 58 |

1. Awani Kumar 2006     D1= 130    D2= 130    D3= 130    D4= 130

| | | | | | | | | | | | | | | | | | | | |
|---|---|---|---|---|---|---|---|---|---|---|---|---|---|---|---|---|---|---|---|
| 1 | 20 | 47 | 62 | | 48 | 61 | 2 | 19 | | 29 | 16 | 51 | 34 | | 52 | 33 | 30 | 15 |
| 28 | 37 | 54 | 11 | | 53 | 12 | 27 | 38 | | 44 | 21 | 6 | 59 | | 5 | 60 | 43 | 22 |
| 45 | 64 | 3 | 18 | | 4 | 17 | 46 | 63 | | 49 | 36 | 31 | 14 | | 32 | 13 | 50 | 35 |
| 56 | 9 | 26 | 39 | | 25 | 40 | 55 | 10 | | 8 | 57 | 42 | 23 | | 41 | 24 | 7 | 58 |

2. Awani Kumar 2006     D1= 102    D2= 166    D3= 94    D4= 158

| | | | | | | | | | | | | | | | | | | | |
|---|---|---|---|---|---|---|---|---|---|---|---|---|---|---|---|---|---|---|---|
| 1 | 24 | 47 | 58 | | 46 | 59 | 4 | 21 | | 23 | 34 | 57 | 16 | | 60 | 13 | 22 | 35 |
| 42 | 63 | 8 | 17 | | 5 | 20 | 43 | 62 | | 64 | 9 | 18 | 39 | | 19 | 38 | 61 | 12 |
| 55 | 2 | 25 | 48 | | 28 | 45 | 54 | 3 | | 33 | 56 | 15 | 26 | | 14 | 27 | 36 | 53 |
| 32 | 41 | 50 | 7 | | 51 | 6 | 29 | 44 | | 10 | 31 | 40 | 49 | | 37 | 52 | 11 | 30 |

3. Guenter Stertenbrink 2003    D1= 66    D2= 194    D3= 130    D4= 130

| | | | | | | | | | | | | | | | | | | | |
|---|---|---|---|---|---|---|---|---|---|---|---|---|---|---|---|---|---|---|---|
| 1 | 42 | 55 | 32 | | 56 | 31 | 34 | 9 | | 47 | 8 | 25 | 50 | | 26 | 49 | 16 | 39 |
| 46 | 5 | 28 | 51 | | 27 | 52 | 13 | 38 | | 4 | 43 | 54 | 29 | | 53 | 30 | 35 | 12 |
| 23 | 64 | 33 | 10 | | 2 | 41 | 24 | 63 | | 57 | 18 | 15 | 40 | | 48 | 7 | 58 | 17 |
| 60 | 19 | 14 | 37 | | 45 | 6 | 59 | 20 | | 22 | 61 | 36 | 11 | | 3 | 44 | 21 | 62 |

4. Guenter Stertenbrink 2003    D1= 130    D2= 66    D3= 194    D4= 130

| | | | | | | | | | | | | | | | | | | | |
|---|---|---|---|---|---|---|---|---|---|---|---|---|---|---|---|---|---|---|---|
| 2 | 11 | 62 | 55 | | 63 | 54 | 3 | 10 | | 14 | 7 | 50 | 59 | | 51 | 58 | 15 | 6 |
| 45 | 56 | 1 | 28 | | 4 | 25 | 48 | 53 | | 49 | 44 | 29 | 8 | | 32 | 5 | 52 | 41 |
| 64 | 37 | 20 | 9 | | 17 | 12 | 61 | 40 | | 36 | 57 | 16 | 21 | | 13 | 24 | 33 | 60 |
| 19 | 26 | 47 | 38 | | 46 | 39 | 18 | 27 | | 31 | 22 | 35 | 42 | | 34 | 43 | 30 | 23 |

5. Awani Kumar 2009     D1= 66    D2= 194    D3= 66    D4= 194

| | | | | | | | | | | | | | | | | | | | |
|---|---|---|---|---|---|---|---|---|---|---|---|---|---|---|---|---|---|---|---|
| 2 | 11 | 62 | 55 | | 63 | 54 | 3 | 10 | | 14 | 23 | 50 | 43 | | 51 | 42 | 15 | 22 |
| 45 | 56 | 1 | 28 | | 4 | 25 | 48 | 53 | | 49 | 44 | 29 | 8 | | 32 | 5 | 52 | 41 |
| 64 | 37 | 20 | 9 | | 17 | 12 | 61 | 40 | | 36 | 57 | 16 | 21 | | 13 | 24 | 33 | 60 |
| 19 | 26 | 47 | 38 | | 46 | 39 | 18 | 27 | | 31 | 6 | 35 | 58 | | 34 | 59 | 30 | 7 |

6. Francis Gaspalou 2009    D1= 50    D2= 194    D3= 82    D4= 194

| | | | | | | | | | | | | | | | | | | | |
|---|---|---|---|---|---|---|---|---|---|---|---|---|---|---|---|---|---|---|---|
| 2 | 11 | 62 | 55 | | 63 | 54 | 3 | 10 | | 30 | 7 | 34 | 59 | | 35 | 58 | 31 | 6 |
| 45 | 56 | 1 | 28 | | 4 | 25 | 48 | 53 | | 49 | 44 | 29 | 8 | | 32 | 5 | 52 | 41 |
| 64 | 37 | 20 | 9 | | 17 | 12 | 61 | 40 | | 36 | 57 | 16 | 21 | | 13 | 24 | 33 | 60 |
| 19 | 26 | 47 | 38 | | 46 | 39 | 18 | 27 | | 15 | 22 | 51 | 42 | | 50 | 43 | 14 | 23 |

7. Francis Gaspalou 2009    D1= 66    D2= 210    D3= 66    D4= 178

| | | | | | | | | | | | | | | | | | | | |
|---|---|---|---|---|---|---|---|---|---|---|---|---|---|---|---|---|---|---|---|
| 2 | 11 | 62 | 55 | | 63 | 54 | 3 | 10 | | 30 | 23 | 34 | 43 | | 35 | 42 | 31 | 22 |
| 45 | 56 | 1 | 28 | | 4 | 25 | 48 | 53 | | 49 | 44 | 29 | 8 | | 32 | 5 | 52 | 41 |
| 64 | 37 | 20 | 9 | | 17 | 12 | 61 | 40 | | 36 | 57 | 16 | 21 | | 13 | 24 | 33 | 60 |
| 19 | 26 | 47 | 38 | | 46 | 39 | 18 | 27 | | 15 | 6 | 51 | 58 | | 50 | 59 | 14 | 7 |

8. Francis Gaspalou 2009    D1= 50    D2= 210    D3= 82    D4= 178

### 9. Awani Kumar 2009

| 2 | 11 | 64 | 53 | | 61 | 56 | 3 | 10 | | 12 | 1 | 54 | 63 | | 55 | 62 | 9 | 4 |
|---|---|---|---|---|---|---|---|---|---|---|---|---|---|---|---|---|---|---|
| 21 | 32 | 43 | 34 | | 42 | 35 | 24 | 29 | | 31 | 22 | 33 | 44 | | 36 | 41 | 30 | 23 |
| 60 | 49 | 6 | 15 | | 7 | 14 | 57 | 52 | | 50 | 59 | 16 | 5 | | 13 | 8 | 51 | 58 |
| 47 | 38 | 17 | 28 | | 20 | 25 | 46 | 39 | | 37 | 48 | 27 | 18 | | 26 | 19 | 40 | 45 |

D1= 98   D2= 162   D3= 98   D4= 162

### 10. Awani Kumar 2009

| 2 | 11 | 64 | 53 | | 61 | 56 | 3 | 10 | | 12 | 1 | 54 | 63 | | 55 | 62 | 9 | 4 |
|---|---|---|---|---|---|---|---|---|---|---|---|---|---|---|---|---|---|---|
| 47 | 38 | 17 | 28 | | 20 | 25 | 46 | 39 | | 37 | 48 | 27 | 18 | | 26 | 19 | 40 | 45 |
| 60 | 49 | 6 | 15 | | 7 | 14 | 57 | 52 | | 50 | 59 | 16 | 5 | | 13 | 8 | 51 | 58 |
| 21 | 32 | 43 | 34 | | 42 | 35 | 24 | 29 | | 31 | 22 | 33 | 44 | | 36 | 41 | 30 | 23 |

D1= 66   D2= 194   D3= 66   D4= 194

### 11. Guenter Stertenbrink 2003

| 2 | 21 | 60 | 47 | | 61 | 42 | 7 | 20 | | 28 | 15 | 34 | 53 | | 39 | 52 | 29 | 10 |
|---|---|---|---|---|---|---|---|---|---|---|---|---|---|---|---|---|---|---|
| 27 | 16 | 33 | 54 | | 40 | 51 | 30 | 9 | | 1 | 22 | 59 | 48 | | 62 | 41 | 8 | 19 |
| 64 | 43 | 6 | 17 | | 3 | 24 | 57 | 46 | | 38 | 49 | 32 | 11 | | 25 | 14 | 35 | 56 |
| 37 | 50 | 31 | 12 | | 26 | 13 | 36 | 55 | | 63 | 44 | 5 | 18 | | 4 | 23 | 58 | 45 |

D1= 130   D2= 130   D3= 130   D4= 130

### 12. Awani Kumar 2017

| 2 | 21 | 60 | 47 | | 61 | 42 | 7 | 20 | | 28 | 15 | 34 | 53 | | 39 | 52 | 29 | 10 |
|---|---|---|---|---|---|---|---|---|---|---|---|---|---|---|---|---|---|---|
| 59 | 48 | 1 | 22 | | 8 | 19 | 62 | 41 | | 33 | 54 | 27 | 16 | | 30 | 9 | 40 | 51 |
| 32 | 43 | 6 | 49 | | 3 | 56 | 25 | 46 | | 38 | 17 | 64 | 11 | | 57 | 14 | 35 | 24 |
| 37 | 18 | 63 | 12 | | 58 | 13 | 36 | 23 | | 31 | 44 | 5 | 50 | | 4 | 55 | 26 | 45 |

D1= 130   D2= 130   D3= 130   D4= 130

### 13. Francis Gaspalou 2009

| 2 | 21 | 60 | 47 | | 61 | 42 | 7 | 20 | | 12 | 31 | 50 | 37 | | 55 | 36 | 13 | 26 |
|---|---|---|---|---|---|---|---|---|---|---|---|---|---|---|---|---|---|---|
| 59 | 48 | 1 | 22 | | 8 | 19 | 62 | 41 | | 49 | 38 | 11 | 32 | | 14 | 25 | 56 | 35 |
| 64 | 43 | 6 | 17 | | 3 | 24 | 57 | 46 | | 54 | 33 | 16 | 27 | | 9 | 30 | 51 | 40 |
| 5 | 18 | 63 | 44 | | 58 | 45 | 4 | 23 | | 15 | 28 | 53 | 34 | | 52 | 39 | 10 | 29 |

D1= 66   D2= 194   D3= 66   D4= 194

### 14. Francis Gaspalou 2009

| 2 | 21 | 60 | 47 | | 61 | 42 | 7 | 20 | | 28 | 15 | 34 | 53 | | 39 | 52 | 29 | 10 |
|---|---|---|---|---|---|---|---|---|---|---|---|---|---|---|---|---|---|---|
| 59 | 48 | 1 | 22 | | 8 | 19 | 62 | 41 | | 33 | 54 | 27 | 16 | | 30 | 9 | 40 | 51 |
| 64 | 43 | 6 | 17 | | 3 | 24 | 57 | 46 | | 38 | 49 | 32 | 11 | | 25 | 14 | 35 | 56 |
| 5 | 18 | 63 | 44 | | 58 | 45 | 4 | 23 | | 31 | 12 | 37 | 50 | | 36 | 55 | 26 | 13 |

D1= 66   D2= 194   D3= 66   D4= 194

### 15. Francis Gaspalou 2009

| 2 | 23 | 64 | 41 | | 61 | 44 | 3 | 22 | | 56 | 33 | 10 | 31 | | 11 | 30 | 53 | 36 |
|---|---|---|---|---|---|---|---|---|---|---|---|---|---|---|---|---|---|---|
| 57 | 48 | 7 | 18 | | 6 | 19 | 60 | 45 | | 15 | 26 | 49 | 40 | | 52 | 37 | 14 | 27 |
| 24 | 1 | 42 | 63 | | 43 | 62 | 21 | 4 | | 34 | 55 | 32 | 9 | | 29 | 12 | 35 | 54 |
| 47 | 58 | 17 | 8 | | 20 | 5 | 46 | 59 | | 25 | 16 | 39 | 50 | | 38 | 51 | 28 | 13 |

D1= 66   D2= 194   D3= 194   D4= 66

### 16. Awani Kumar 2017

| 2 | 25 | 56 | 47 | | 61 | 38 | 11 | 20 | | 24 | 15 | 34 | 57 | | 43 | 52 | 29 | 6 |
|---|---|---|---|---|---|---|---|---|---|---|---|---|---|---|---|---|---|---|
| 55 | 48 | 1 | 26 | | 12 | 19 | 62 | 37 | | 33 | 58 | 23 | 16 | | 30 | 5 | 44 | 51 |
| 32 | 39 | 10 | 49 | | 3 | 60 | 21 | 46 | | 42 | 17 | 64 | 7 | | 53 | 14 | 35 | 28 |
| 41 | 18 | 63 | 8 | | 54 | 13 | 36 | 27 | | 31 | 40 | 9 | 50 | | 4 | 59 | 22 | 45 |

D1= 66   D2= 130   D3= 130   D4= 130

### 17. Francis Gaspalou 2009

| 2 | 27 | 46 | 55 | | 47 | 54 | 3 | 26 | | 62 | 39 | 18 | 11 | | 19 | 10 | 63 | 38 |
|---|---|---|---|---|---|---|---|---|---|---|---|---|---|---|---|---|---|---|
| 29 | 40 | 49 | 12 | | 52 | 9 | 32 | 37 | | 33 | 28 | 13 | 56 | | 16 | 53 | 36 | 25 |
| 48 | 21 | 4 | 57 | | 1 | 60 | 45 | 24 | | 20 | 41 | 64 | 5 | | 61 | 8 | 17 | 44 |
| 51 | 42 | 31 | 6 | | 30 | 7 | 50 | 43 | | 15 | 22 | 35 | 58 | | 34 | 59 | 14 | 23 |

D1= 98   D2= 162   D3= 162   D4= 98

### 18. Awani Kumar 2007

| 2  | 27 | 46 | 55 |   | 47 | 54 | 3  | 26 |   | 30 | 7  | 50 | 43 |   | 51 | 42 | 31 | 6  |
|----|----|----|----|---|----|----|----|----|---|----|----|----|----|---|----|----|----|----|
| 45 | 56 | 1  | 28 |   | 4  | 25 | 48 | 53 |   | 49 | 44 | 29 | 8  |   | 32 | 5  | 52 | 41 |
| 64 | 37 | 20 | 9  |   | 17 | 12 | 61 | 40 |   | 36 | 57 | 16 | 21 |   | 13 | 24 | 33 | 60 |
| 19 | 10 | 63 | 38 |   | 62 | 39 | 18 | 11 |   | 15 | 22 | 35 | 58 |   | 34 | 59 | 14 | 23 |

D1= 66   D2= 194   D3= 66   D4= 194

### 19. Francis Gaspalou 2009

| 2  | 27 | 46 | 55 |   | 61 | 44 | 17 | 8  |   | 16 | 53 | 36 | 25 |   | 51 | 6  | 31 | 42 |
|----|----|----|----|---|----|----|----|----|---|----|----|----|----|---|----|----|----|----|
| 47 | 54 | 3  | 26 |   | 20 | 5  | 64 | 41 |   | 33 | 28 | 13 | 56 |   | 30 | 43 | 50 | 7  |
| 62 | 11 | 18 | 39 |   | 1  | 60 | 45 | 24 |   | 52 | 37 | 32 | 9  |   | 15 | 22 | 35 | 58 |
| 19 | 38 | 63 | 10 |   | 48 | 21 | 4  | 57 |   | 29 | 12 | 49 | 40 |   | 34 | 59 | 14 | 23 |

D1= 62   D2= 190   D3= 134   D4= 134

### 20. Francis Gaspalou 2009

| 2  | 27 | 46 | 55 |   | 61 | 8  | 17 | 44 |   | 16 | 53 | 36 | 25 |   | 51 | 42 | 31 | 6  |
|----|----|----|----|---|----|----|----|----|---|----|----|----|----|---|----|----|----|----|
| 47 | 54 | 3  | 26 |   | 20 | 41 | 64 | 5  |   | 33 | 28 | 13 | 56 |   | 30 | 7  | 50 | 43 |
| 62 | 39 | 18 | 11 |   | 1  | 60 | 45 | 24 |   | 52 | 9  | 32 | 37 |   | 15 | 22 | 35 | 58 |
| 19 | 10 | 63 | 38 |   | 48 | 21 | 4  | 57 |   | 29 | 40 | 49 | 12 |   | 34 | 59 | 14 | 23 |

D1= 98   D2= 162   D3= 98   D4= 162

### 21. Francis Gaspalou 2009

| 2  | 27 | 62 | 39 |   | 63 | 38 | 3  | 26 |   | 14 | 7  | 50 | 59 |   | 51 | 58 | 15 | 6  |
|----|----|----|----|---|----|----|----|----|---|----|----|----|----|---|----|----|----|----|
| 45 | 56 | 1  | 28 |   | 4  | 25 | 48 | 53 |   | 49 | 44 | 29 | 8  |   | 32 | 5  | 52 | 41 |
| 64 | 37 | 20 | 9  |   | 17 | 12 | 61 | 40 |   | 36 | 57 | 16 | 21 |   | 13 | 24 | 33 | 60 |
| 19 | 10 | 47 | 54 |   | 46 | 55 | 18 | 11 |   | 31 | 22 | 35 | 42 |   | 34 | 43 | 30 | 23 |

D1= 66   D2= 178   D3= 66   D4= 210

### 22. Francis Gaspalou 2009

| 2  | 27 | 62 | 39 |   | 63 | 38 | 3  | 26 |   | 14 | 23 | 50 | 43 |   | 51 | 42 | 15 | 22 |
|----|----|----|----|---|----|----|----|----|---|----|----|----|----|---|----|----|----|----|
| 45 | 56 | 1  | 28 |   | 4  | 25 | 48 | 53 |   | 49 | 44 | 29 | 8  |   | 32 | 5  | 52 | 41 |
| 64 | 37 | 20 | 9  |   | 17 | 12 | 61 | 40 |   | 36 | 57 | 16 | 21 |   | 13 | 24 | 33 | 60 |
| 19 | 10 | 47 | 54 |   | 46 | 55 | 18 | 11 |   | 31 | 6  | 35 | 58 |   | 34 | 59 | 30 | 7  |

D1= 50   D2= 178   D3= 82   D4= 210

### 23. Francis Gaspalou 2009

| 2  | 27 | 62 | 39 |   | 63 | 38 | 3  | 26 |   | 30 | 7  | 34 | 59 |   | 35 | 58 | 31 | 6  |
|----|----|----|----|---|----|----|----|----|---|----|----|----|----|---|----|----|----|----|
| 45 | 56 | 1  | 28 |   | 4  | 25 | 48 | 53 |   | 49 | 44 | 29 | 8  |   | 32 | 5  | 52 | 41 |
| 64 | 37 | 20 | 9  |   | 17 | 12 | 61 | 40 |   | 36 | 57 | 16 | 21 |   | 13 | 24 | 33 | 60 |
| 19 | 10 | 47 | 54 |   | 46 | 55 | 18 | 11 |   | 15 | 22 | 51 | 42 |   | 50 | 43 | 14 | 23 |

D1= 66   D2= 194   D3= 66   D4= 194

### 24. Francis Gaspalou 2009

| 2  | 27 | 62 | 39 |   | 63 | 38 | 3  | 26 |   | 30 | 23 | 34 | 43 |   | 35 | 42 | 31 | 22 |
|----|----|----|----|---|----|----|----|----|---|----|----|----|----|---|----|----|----|----|
| 45 | 56 | 1  | 28 |   | 4  | 25 | 48 | 53 |   | 49 | 44 | 29 | 8  |   | 32 | 5  | 52 | 41 |
| 64 | 37 | 20 | 9  |   | 17 | 12 | 61 | 40 |   | 36 | 57 | 16 | 21 |   | 13 | 24 | 33 | 60 |
| 19 | 10 | 47 | 54 |   | 46 | 55 | 18 | 11 |   | 15 | 6  | 51 | 58 |   | 50 | 59 | 14 | 7  |

D1= 50   D2= 194   D3= 82   D4= 194

### 25. Awani Kumar 2017

| 2  | 27 | 64 | 37 |   | 61 | 40 | 19 | 10 |   | 28 | 33 | 38 | 31 |   | 39 | 30 | 9  | 52 |
|----|----|----|----|---|----|----|----|----|---|----|----|----|----|---|----|----|----|----|
| 45 | 56 | 3  | 26 |   | 18 | 11 | 48 | 53 |   | 55 | 14 | 25 | 36 |   | 12 | 49 | 54 | 15 |
| 60 | 1  | 6  | 63 |   | 7  | 62 | 41 | 20 |   | 34 | 59 | 32 | 5  |   | 29 | 8  | 51 | 42 |
| 23 | 46 | 57 | 4  |   | 44 | 17 | 22 | 47 |   | 13 | 24 | 35 | 58 |   | 50 | 43 | 16 | 21 |

D1= 66   D2= 194   D3= 162   D4= 98

### 26. Francis Gaspalou 2009

| 2  | 29 | 48 | 51 |   | 59 | 40 | 21 | 10 |   | 30 | 49 | 4  | 47 |   | 39 | 12 | 57 | 22 |
|----|----|----|----|---|----|----|----|----|---|----|----|----|----|---|----|----|----|----|
| 31 | 52 | 1  | 46 |   | 38 | 9  | 60 | 23 |   | 3  | 32 | 45 | 50 |   | 58 | 37 | 24 | 11 |
| 62 | 33 | 20 | 15 |   | 7  | 28 | 41 | 54 |   | 34 | 13 | 64 | 19 |   | 27 | 56 | 5  | 42 |
| 35 | 16 | 61 | 18 |   | 26 | 53 | 8  | 43 |   | 63 | 36 | 17 | 14 |   | 6  | 25 | 44 | 55 |

D1= 130   D2= 130   D3= 130   D4= 130

| 2 | 29 | 48 | 51 | | 55 | 44 | 25 | 6 | | 12 | 49 | 38 | 31 | | 61 | 8 | 19 | 42 |
|---|---|---|---|---|---|---|---|---|---|---|---|---|---|---|---|---|---|---|
| 47 | 52 | 1 | 30 | | 26 | 5 | 56 | 43 | | 37 | 32 | 11 | 50 | | 20 | 41 | 62 | 7 |
| 54 | 33 | 28 | 15 | | 3 | 24 | 45 | 58 | | 64 | 13 | 18 | 35 | | 9 | 60 | 39 | 22 |
| 27 | 16 | 53 | 34 | | 46 | 57 | 4 | 23 | | 17 | 36 | 63 | 14 | | 40 | 21 | 10 | 59 |

27. Awani Kumar 2012    D1= 84    D2= 160    D3= 104    D4= 172

| 2 | 29 | 48 | 51 | | 59 | 12 | 21 | 38 | | 46 | 49 | 4 | 31 | | 23 | 40 | 57 | 10 |
|---|---|---|---|---|---|---|---|---|---|---|---|---|---|---|---|---|---|---|
| 47 | 52 | 1 | 30 | | 22 | 37 | 60 | 11 | | 3 | 32 | 45 | 50 | | 58 | 9 | 24 | 39 |
| 62 | 33 | 20 | 15 | | 43 | 28 | 5 | 54 | | 18 | 13 | 64 | 35 | | 7 | 56 | 41 | 26 |
| 19 | 16 | 61 | 34 | | 6 | 53 | 44 | 27 | | 63 | 36 | 17 | 14 | | 42 | 25 | 8 | 55 |

28. Awani Kumar 2007    D1= 158    D2= 166    D3= 102    D4= 94

| 2 | 29 | 48 | 51 | | 59 | 40 | 21 | 10 | | 46 | 49 | 4 | 31 | | 23 | 12 | 57 | 38 |
|---|---|---|---|---|---|---|---|---|---|---|---|---|---|---|---|---|---|---|
| 47 | 52 | 1 | 30 | | 22 | 9 | 60 | 39 | | 3 | 32 | 45 | 50 | | 58 | 37 | 24 | 11 |
| 62 | 33 | 20 | 15 | | 7 | 28 | 41 | 54 | | 18 | 13 | 64 | 35 | | 43 | 56 | 5 | 26 |
| 19 | 16 | 61 | 34 | | 42 | 53 | 8 | 27 | | 63 | 36 | 17 | 14 | | 6 | 25 | 44 | 55 |

29. Awani Kumar 2007    D1= 130    D2= 130    D3= 130    D4= 130

| 2 | 29 | 56 | 43 | | 57 | 38 | 15 | 20 | | 24 | 11 | 34 | 61 | | 47 | 52 | 25 | 6 |
|---|---|---|---|---|---|---|---|---|---|---|---|---|---|---|---|---|---|---|
| 55 | 44 | 1 | 30 | | 16 | 19 | 58 | 37 | | 33 | 62 | 23 | 12 | | 26 | 5 | 48 | 51 |
| 64 | 3 | 42 | 21 | | 39 | 28 | 49 | 14 | | 10 | 53 | 32 | 35 | | 17 | 46 | 7 | 60 |
| 9 | 54 | 31 | 36 | | 18 | 45 | 8 | 59 | | 63 | 4 | 41 | 22 | | 40 | 27 | 50 | 13 |

30. Awani Kumar 2017    D1= 66    D2= 194    D3= 66    D4= 194

| 2 | 29 | 60 | 39 | | 61 | 34 | 7 | 28 | | 52 | 47 | 10 | 21 | | 15 | 20 | 53 | 42 |
|---|---|---|---|---|---|---|---|---|---|---|---|---|---|---|---|---|---|---|
| 51 | 48 | 1 | 30 | | 16 | 19 | 62 | 33 | | 9 | 22 | 59 | 40 | | 54 | 41 | 8 | 27 |
| 64 | 3 | 6 | 57 | | 35 | 32 | 25 | 38 | | 14 | 49 | 56 | 11 | | 17 | 46 | 43 | 24 |
| 13 | 50 | 63 | 4 | | 18 | 45 | 36 | 31 | | 55 | 12 | 5 | 58 | | 44 | 23 | 26 | 37 |

31. Awani Kumar 2017    D1= 114    D2= 194    D3= 146    D4= 66

| 2 | 29 | 60 | 39 | | 59 | 40 | 1 | 30 | | 32 | 3 | 38 | 57 | | 37 | 58 | 31 | 4 |
|---|---|---|---|---|---|---|---|---|---|---|---|---|---|---|---|---|---|---|
| 53 | 42 | 15 | 20 | | 16 | 19 | 54 | 41 | | 43 | 56 | 17 | 14 | | 18 | 13 | 44 | 55 |
| 28 | 7 | 34 | 61 | | 33 | 62 | 27 | 8 | | 6 | 25 | 64 | 35 | | 63 | 36 | 5 | 26 |
| 47 | 52 | 21 | 10 | | 22 | 9 | 48 | 51 | | 49 | 46 | 11 | 24 | | 12 | 23 | 50 | 45 |

32. Guenter Stertenbrink 2003    D1= 130    D2= 130    D3= 130    D4= 130

| 2 | 29 | 60 | 39 | | 59 | 40 | 1 | 30 | | 64 | 3 | 38 | 25 | | 5 | 58 | 31 | 36 |
|---|---|---|---|---|---|---|---|---|---|---|---|---|---|---|---|---|---|---|
| 53 | 42 | 15 | 20 | | 16 | 19 | 54 | 41 | | 43 | 24 | 49 | 14 | | 18 | 45 | 12 | 55 |
| 28 | 7 | 34 | 61 | | 33 | 62 | 27 | 8 | | 6 | 57 | 32 | 35 | | 63 | 4 | 37 | 26 |
| 47 | 52 | 21 | 10 | | 22 | 9 | 48 | 51 | | 17 | 46 | 11 | 56 | | 44 | 23 | 50 | 13 |

33. Awani Kumar 2017    D1= 66    D2= 194    D3= 194    D4= 66

| 2 | 29 | 60 | 39 | | 61 | 34 | 7 | 28 | | 52 | 47 | 10 | 21 | | 15 | 20 | 53 | 42 |
|---|---|---|---|---|---|---|---|---|---|---|---|---|---|---|---|---|---|---|
| 59 | 40 | 9 | 22 | | 8 | 27 | 54 | 41 | | 1 | 30 | 51 | 48 | | 62 | 33 | 16 | 19 |
| 64 | 3 | 6 | 57 | | 35 | 32 | 25 | 38 | | 14 | 49 | 56 | 11 | | 17 | 46 | 43 | 24 |
| 5 | 58 | 55 | 12 | | 26 | 37 | 44 | 23 | | 63 | 4 | 13 | 50 | | 36 | 31 | 18 | 45 |

34. Awani Kumar 2017    D1= 130    D2= 178    D3= 130    D4= 82

| 2 | 31 | 54 | 43 | | 61 | 36 | 9 | 24 | | 8 | 17 | 60 | 45 | | 59 | 46 | 7 | 18 |
|---|---|---|---|---|---|---|---|---|---|---|---|---|---|---|---|---|---|---|
| 55 | 42 | 3 | 30 | | 12 | 21 | 64 | 33 | | 57 | 48 | 5 | 20 | | 6 | 19 | 58 | 47 |
| 62 | 35 | 10 | 23 | | 1 | 32 | 53 | 44 | | 52 | 37 | 16 | 25 | | 15 | 26 | 51 | 38 |
| 11 | 22 | 63 | 34 | | 56 | 41 | 4 | 29 | | 13 | 28 | 49 | 40 | | 50 | 39 | 14 | 27 |

35. Awani Kumar 2011    D1= 66    D2= 194    D3= 66    D4= 194

| 2 | 31 | 60 | 37 |
|---|---|---|---|
| 59 | 38 | 1 | 32 |
| 64 | 33 | 6 | 27 |
| 5 | 28 | 63 | 34 |

| 61 | 36 | 7 | 26 |
|---|---|---|---|
| 8 | 25 | 62 | 35 |
| 3 | 30 | 57 | 40 |
| 58 | 39 | 4 | 29 |

| 12 | 53 | 50 | 15 |
|---|---|---|---|
| 17 | 48 | 43 | 22 |
| 54 | 11 | 16 | 49 |
| 47 | 18 | 21 | 44 |

| 55 | 10 | 13 | 52 |
|---|---|---|---|
| 46 | 19 | 24 | 41 |
| 9 | 56 | 51 | 14 |
| 20 | 45 | 42 | 23 |

36. Awani Kumar 2017   D1= 66   D2= 130   D3= 130   D4= 194

| 2 | 39 | 60 | 29 |
|---|---|---|---|
| 57 | 32 | 7 | 34 |
| 40 | 17 | 30 | 43 |
| 31 | 42 | 33 | 24 |

| 61 | 28 | 3 | 38 |
|---|---|---|---|
| 6 | 35 | 64 | 25 |
| 27 | 46 | 37 | 20 |
| 36 | 21 | 26 | 47 |

| 56 | 1 | 14 | 59 |
|---|---|---|---|
| 15 | 58 | 49 | 8 |
| 18 | 55 | 44 | 13 |
| 41 | 16 | 23 | 50 |

| 11 | 62 | 53 | 4 |
|---|---|---|---|
| 52 | 5 | 10 | 63 |
| 45 | 12 | 19 | 54 |
| 22 | 51 | 48 | 9 |

37. Awani Kumar 2012   D1= 90   D2= 170   D3= 130   D4= 130

| 2 | 43 | 54 | 31 |
|---|---|---|---|
| 55 | 30 | 3 | 42 |
| 62 | 23 | 10 | 35 |
| 11 | 34 | 63 | 22 |

| 61 | 24 | 9 | 36 |
|---|---|---|---|
| 12 | 33 | 64 | 21 |
| 1 | 44 | 53 | 32 |
| 56 | 29 | 4 | 41 |

| 8 | 45 | 60 | 17 |
|---|---|---|---|
| 57 | 20 | 5 | 48 |
| 52 | 25 | 16 | 37 |
| 13 | 40 | 49 | 28 |

| 59 | 18 | 7 | 46 |
|---|---|---|---|
| 6 | 47 | 58 | 19 |
| 15 | 38 | 51 | 26 |
| 50 | 27 | 14 | 39 |

38. Awani Kumar 2012   D1= 90   D2= 170   D3= 106   D4= 154

| 2 | 43 | 54 | 31 |
|---|---|---|---|
| 55 | 30 | 3 | 42 |
| 62 | 23 | 10 | 35 |
| 11 | 34 | 63 | 22 |

| 61 | 24 | 9 | 36 |
|---|---|---|---|
| 12 | 33 | 64 | 21 |
| 1 | 44 | 53 | 32 |
| 56 | 29 | 4 | 41 |

| 16 | 37 | 52 | 25 |
|---|---|---|---|
| 49 | 28 | 13 | 40 |
| 60 | 17 | 8 | 45 |
| 5 | 48 | 57 | 20 |

| 51 | 26 | 15 | 38 |
|---|---|---|---|
| 14 | 39 | 50 | 27 |
| 7 | 46 | 59 | 18 |
| 58 | 19 | 6 | 47 |

39. Awani Kumar 2012   D1= 90   D2= 170   D3= 106   D4= 154

| 2 | 43 | 62 | 23 |
|---|---|---|---|
| 61 | 8 | 17 | 44 |
| 16 | 53 | 36 | 25 |
| 51 | 26 | 15 | 38 |

| 63 | 22 | 3 | 42 |
|---|---|---|---|
| 20 | 41 | 64 | 5 |
| 33 | 28 | 13 | 56 |
| 14 | 39 | 50 | 27 |

| 46 | 7 | 18 | 59 |
|---|---|---|---|
| 1 | 60 | 45 | 24 |
| 52 | 9 | 32 | 37 |
| 31 | 54 | 35 | 10 |

| 19 | 58 | 47 | 6 |
|---|---|---|---|
| 48 | 21 | 4 | 57 |
| 29 | 40 | 49 | 12 |
| 34 | 11 | 30 | 55 |

40. Awani Kumar 2009   D1= 130   D2= 130   D3= 130   D4= 130

| 2 | 43 | 62 | 23 |
|---|---|---|---|
| 61 | 24 | 9 | 36 |
| 16 | 37 | 52 | 25 |
| 51 | 26 | 7 | 46 |

| 63 | 22 | 3 | 42 |
|---|---|---|---|
| 12 | 33 | 64 | 21 |
| 49 | 28 | 13 | 40 |
| 6 | 47 | 50 | 27 |

| 54 | 31 | 10 | 35 |
|---|---|---|---|
| 1 | 44 | 53 | 32 |
| 60 | 17 | 8 | 45 |
| 15 | 38 | 59 | 18 |

| 11 | 34 | 55 | 30 |
|---|---|---|---|
| 56 | 29 | 4 | 41 |
| 5 | 48 | 57 | 20 |
| 58 | 19 | 14 | 39 |

41. Awani Kumar 2012   D1= 82   D2= 162   D3= 162   D4= 114

| 2 | 43 | 62 | 23 |
|---|---|---|---|
| 61 | 24 | 1 | 44 |
| 16 | 53 | 36 | 25 |
| 51 | 10 | 31 | 38 |

| 63 | 22 | 3 | 42 |
|---|---|---|---|
| 4 | 41 | 64 | 21 |
| 33 | 28 | 13 | 56 |
| 30 | 39 | 50 | 11 |

| 46 | 7 | 18 | 59 |
|---|---|---|---|
| 17 | 60 | 45 | 8 |
| 52 | 9 | 32 | 37 |
| 15 | 54 | 35 | 26 |

| 19 | 58 | 47 | 6 |
|---|---|---|---|
| 48 | 5 | 20 | 57 |
| 29 | 40 | 49 | 12 |
| 34 | 27 | 14 | 55 |

42. Awani Kumar 2009   D1= 130   D2= 130   D3= 130   D4= 130

| 2 | 43 | 64 | 21 |
|---|---|---|---|
| 61 | 24 | 3 | 42 |
| 28 | 33 | 38 | 31 |
| 39 | 30 | 25 | 36 |

| 63 | 22 | 1 | 44 |
|---|---|---|---|
| 4 | 41 | 62 | 23 |
| 53 | 16 | 11 | 50 |
| 10 | 51 | 56 | 13 |

| 60 | 17 | 6 | 47 |
|---|---|---|---|
| 7 | 46 | 57 | 20 |
| 34 | 27 | 32 | 37 |
| 29 | 40 | 35 | 26 |

| 5 | 48 | 59 | 18 |
|---|---|---|---|
| 58 | 19 | 8 | 45 |
| 15 | 54 | 49 | 12 |
| 52 | 9 | 14 | 55 |

43. Awani Kumar 2012   D1= 130   D2= 162   D3= 130   D4= 98

| 2 | 45 | 60 | 23 |
|---|---|---|---|
| 59 | 24 | 33 | 14 |
| 64 | 3 | 6 | 57 |
| 5 | 58 | 31 | 36 |

| 61 | 18 | 7 | 44 |
|---|---|---|---|
| 8 | 43 | 30 | 49 |
| 19 | 48 | 41 | 22 |
| 42 | 21 | 52 | 15 |

| 28 | 55 | 34 | 13 |
|---|---|---|---|
| 1 | 46 | 27 | 56 |
| 38 | 25 | 32 | 35 |
| 63 | 4 | 37 | 26 |

| 39 | 12 | 29 | 50 |
|---|---|---|---|
| 62 | 17 | 40 | 11 |
| 9 | 54 | 51 | 16 |
| 20 | 47 | 10 | 53 |

44. Awani Kumar 2017   D1= 130   D2= 98   D3= 130   D4= 162

| | | | | | | | | | | | | | | | |
|---|---|---|---|---|---|---|---|---|---|---|---|---|---|---|---|
| 2 | 45 | 64 | 19 | 55 | 28 | 37 | 10 | 46 | 1 | 20 | 63 | 27 | 56 | 9 | 38 |
| 47 | 4 | 17 | 62 | 26 | 53 | 12 | 39 | 3 | 48 | 61 | 18 | 54 | 25 | 40 | 11 |
| 30 | 49 | 36 | 15 | 7 | 44 | 21 | 58 | 50 | 29 | 16 | 35 | 43 | 8 | 57 | 22 |
| 51 | 32 | 13 | 34 | 42 | 5 | 60 | 23 | 31 | 52 | 33 | 14 | 6 | 41 | 24 | 59 |

45. Awani Kumar 2007    D1= 130    D2= 66    D3= 194    D4= 130

| | | | | | | | | | | | | | | | |
|---|---|---|---|---|---|---|---|---|---|---|---|---|---|---|---|
| 2 | 47 | 62 | 19 | 59 | 22 | 7 | 42 | 46 | 3 | 18 | 63 | 23 | 58 | 43 | 6 |
| 51 | 30 | 15 | 34 | 10 | 39 | 54 | 27 | 31 | 50 | 35 | 14 | 38 | 11 | 26 | 55 |
| 48 | 1 | 20 | 61 | 21 | 60 | 41 | 8 | 4 | 45 | 64 | 17 | 57 | 24 | 5 | 44 |
| 29 | 52 | 33 | 16 | 40 | 9 | 28 | 53 | 49 | 32 | 13 | 36 | 12 | 37 | 56 | 25 |

46. Awani Kumar 2006    D1= 130    D2= 130    D3= 130    D4= 130

| | | | | | | | | | | | | | | | |
|---|---|---|---|---|---|---|---|---|---|---|---|---|---|---|---|
| 2 | 47 | 62 | 19 | 59 | 22 | 43 | 6 | 46 | 3 | 18 | 63 | 23 | 58 | 7 | 42 |
| 51 | 30 | 15 | 34 | 38 | 11 | 54 | 27 | 31 | 50 | 35 | 14 | 10 | 39 | 26 | 55 |
| 48 | 1 | 20 | 61 | 21 | 60 | 5 | 44 | 4 | 45 | 64 | 17 | 57 | 24 | 41 | 8 |
| 29 | 52 | 33 | 16 | 12 | 37 | 28 | 53 | 49 | 32 | 13 | 36 | 40 | 9 | 56 | 25 |

47. Awani Kumar 2006    D1= 102    D2= 158    D3= 166    D4= 94

| | | | | | | | | | | | | | | | |
|---|---|---|---|---|---|---|---|---|---|---|---|---|---|---|---|
| 2 | 47 | 62 | 19 | 61 | 20 | 1 | 48 | 16 | 33 | 52 | 29 | 51 | 30 | 15 | 34 |
| 59 | 22 | 7 | 42 | 8 | 41 | 60 | 21 | 53 | 28 | 9 | 40 | 10 | 39 | 54 | 27 |
| 46 | 3 | 18 | 63 | 17 | 64 | 45 | 4 | 36 | 13 | 32 | 49 | 31 | 50 | 35 | 14 |
| 23 | 58 | 43 | 6 | 44 | 5 | 24 | 57 | 25 | 56 | 37 | 12 | 38 | 11 | 26 | 55 |

48. Francis Gaspalou 2009    D1= 130    D2= 130    D3= 130    D4= 130

| | | | | | | | | | | | | | | | |
|---|---|---|---|---|---|---|---|---|---|---|---|---|---|---|---|
| 2 | 55 | 60 | 13 | 61 | 12 | 3 | 54 | 56 | 17 | 14 | 43 | 11 | 46 | 53 | 20 |
| 57 | 16 | 7 | 50 | 6 | 51 | 64 | 9 | 15 | 42 | 49 | 24 | 52 | 21 | 10 | 47 |
| 40 | 1 | 30 | 59 | 27 | 62 | 37 | 4 | 18 | 39 | 44 | 29 | 45 | 28 | 19 | 38 |
| 31 | 58 | 33 | 8 | 36 | 5 | 26 | 63 | 41 | 32 | 23 | 34 | 22 | 35 | 48 | 25 |

49. Awani Kumar 2012    D1= 122    D2= 138    D3= 162    D4= 98

| | | | | | | | | | | | | | | | |
|---|---|---|---|---|---|---|---|---|---|---|---|---|---|---|---|
| 2 | 55 | 64 | 9 | 61 | 12 | 3 | 54 | 56 | 1 | 10 | 63 | 11 | 62 | 53 | 4 |
| 57 | 16 | 7 | 50 | 6 | 51 | 60 | 13 | 15 | 58 | 49 | 8 | 52 | 5 | 14 | 59 |
| 24 | 33 | 42 | 31 | 43 | 30 | 21 | 36 | 34 | 23 | 32 | 41 | 29 | 44 | 35 | 22 |
| 47 | 26 | 17 | 40 | 20 | 37 | 46 | 27 | 25 | 48 | 39 | 18 | 38 | 19 | 28 | 45 |

50. Francis Gaspalou 2009    D1= 130    D2= 130    D3= 130    D4= 130

| | | | | | | | | | | | | | | | |
|---|---|---|---|---|---|---|---|---|---|---|---|---|---|---|---|
| 2 | 59 | 62 | 7 | 63 | 6 | 3 | 58 | 46 | 23 | 18 | 43 | 19 | 42 | 47 | 22 |
| 61 | 8 | 17 | 44 | 20 | 41 | 64 | 5 | 1 | 60 | 45 | 24 | 48 | 21 | 4 | 57 |
| 16 | 53 | 36 | 25 | 33 | 28 | 13 | 56 | 52 | 9 | 32 | 37 | 29 | 40 | 49 | 12 |
| 51 | 10 | 15 | 54 | 14 | 55 | 50 | 11 | 31 | 38 | 35 | 26 | 34 | 27 | 30 | 39 |

51. Awani Kumar 2009    D1= 114    D2= 114    D3= 146    D4= 146

| | | | | | | | | | | | | | | | |
|---|---|---|---|---|---|---|---|---|---|---|---|---|---|---|---|
| 3 | 10 | 63 | 54 | 64 | 37 | 20 | 9 | 45 | 56 | 1 | 28 | 18 | 27 | 46 | 39 |
| 62 | 55 | 2 | 11 | 17 | 12 | 61 | 40 | 4 | 25 | 48 | 53 | 47 | 38 | 19 | 26 |
| 15 | 6 | 51 | 58 | 36 | 57 | 16 | 21 | 49 | 44 | 29 | 8 | 30 | 23 | 34 | 43 |
| 50 | 59 | 14 | 7 | 13 | 24 | 33 | 60 | 32 | 5 | 52 | 41 | 35 | 42 | 31 | 22 |

52. Awani Kumar 2009    D1= 66    D2= 194    D3= 194    D4= 66

| | | | | | | | | | | | | | | | |
|---|---|---|---|---|---|---|---|---|---|---|---|---|---|---|---|
| 3 | 10 | 63 | 54 | 64 | 37 | 20 | 9 | 45 | 56 | 1 | 28 | 18 | 27 | 46 | 39 |
| 62 | 55 | 2 | 11 | 17 | 12 | 61 | 40 | 4 | 25 | 48 | 53 | 47 | 38 | 19 | 26 |
| 15 | 22 | 51 | 42 | 36 | 57 | 16 | 21 | 49 | 44 | 29 | 8 | 30 | 7 | 34 | 59 |
| 50 | 43 | 14 | 23 | 13 | 24 | 33 | 60 | 32 | 5 | 52 | 41 | 35 | 58 | 31 | 6 |

53. Francis Gaspalou 2009    D1= 50    D2= 194    D3= 194    D4= 82

### 54. Francis Gaspalou 2009

| 3 | 10 | 63 | 54 | | 64 | 37 | 20 | 9 | | 45 | 56 | 1 | 28 | | 18 | 27 | 46 | 39 |
|---|---|---|---|---|---|---|---|---|---|---|---|---|---|---|---|---|---|---|
| 62 | 55 | 2 | 11 | | 17 | 12 | 61 | 40 | | 4 | 25 | 48 | 53 | | 47 | 38 | 19 | 26 |
| 31 | 6 | 35 | 58 | | 36 | 57 | 16 | 21 | | 49 | 44 | 29 | 8 | | 14 | 23 | 50 | 43 |
| 34 | 59 | 30 | 7 | | 13 | 24 | 33 | 60 | | 32 | 5 | 52 | 41 | | 51 | 42 | 15 | 22 |

D1= 66  D2= 210  D3= 178  D4= 66

### 55. Francis Gaspalou 2009

| 3 | 10 | 63 | 54 | | 64 | 37 | 20 | 9 | | 45 | 56 | 1 | 28 | | 18 | 27 | 46 | 39 |
|---|---|---|---|---|---|---|---|---|---|---|---|---|---|---|---|---|---|---|
| 62 | 55 | 2 | 11 | | 17 | 12 | 61 | 40 | | 4 | 25 | 48 | 53 | | 47 | 38 | 19 | 26 |
| 31 | 22 | 35 | 42 | | 36 | 57 | 16 | 21 | | 49 | 44 | 29 | 8 | | 14 | 7 | 50 | 59 |
| 34 | 43 | 30 | 23 | | 13 | 24 | 33 | 60 | | 32 | 5 | 52 | 41 | | 51 | 58 | 15 | 6 |

D1= 50  D2= 210  D3= 178  D4= 82

### 56. Awani Kumar 2006

| 3 | 18 | 45 | 64 | | 46 | 63 | 4 | 17 | | 31 | 14 | 49 | 36 | | 50 | 35 | 32 | 13 |
|---|---|---|---|---|---|---|---|---|---|---|---|---|---|---|---|---|---|---|
| 26 | 11 | 56 | 37 | | 55 | 38 | 25 | 12 | | 6 | 23 | 44 | 57 | | 43 | 58 | 5 | 24 |
| 47 | 62 | 1 | 20 | | 2 | 19 | 48 | 61 | | 51 | 34 | 29 | 16 | | 30 | 15 | 52 | 33 |
| 54 | 39 | 28 | 9 | | 27 | 10 | 53 | 40 | | 42 | 59 | 8 | 21 | | 7 | 22 | 41 | 60 |

D1= 130  D2= 130  D3= 130  D4= 130

### 57. Awani Kumar 2006

| 3 | 18 | 45 | 64 | | 54 | 11 | 28 | 37 | | 47 | 62 | 1 | 20 | | 26 | 39 | 56 | 9 |
|---|---|---|---|---|---|---|---|---|---|---|---|---|---|---|---|---|---|---|
| 46 | 63 | 4 | 17 | | 27 | 38 | 53 | 12 | | 2 | 19 | 48 | 61 | | 55 | 10 | 25 | 40 |
| 31 | 14 | 49 | 36 | | 6 | 59 | 44 | 21 | | 51 | 34 | 29 | 16 | | 42 | 23 | 8 | 57 |
| 50 | 35 | 32 | 13 | | 43 | 22 | 5 | 60 | | 30 | 15 | 52 | 33 | | 7 | 58 | 41 | 24 |

D1= 94  D2= 158  D3= 166  D4= 102

### 58. Awani Kumar 2017

| 3 | 22 | 47 | 58 | | 46 | 59 | 2 | 23 | | 63 | 38 | 19 | 10 | | 18 | 11 | 62 | 39 |
|---|---|---|---|---|---|---|---|---|---|---|---|---|---|---|---|---|---|---|
| 32 | 37 | 20 | 41 | | 49 | 12 | 61 | 8 | | 36 | 21 | 48 | 25 | | 13 | 60 | 1 | 56 |
| 45 | 28 | 33 | 24 | | 4 | 53 | 16 | 57 | | 17 | 44 | 29 | 40 | | 64 | 5 | 52 | 9 |
| 50 | 43 | 30 | 7 | | 31 | 6 | 51 | 42 | | 14 | 27 | 34 | 55 | | 35 | 54 | 15 | 26 |

D1= 70  D2= 198  D3= 190  D4= 62

### 59. Francis Gaspalou 2009

| 3 | 22 | 47 | 58 | | 64 | 37 | 20 | 9 | | 45 | 28 | 1 | 56 | | 18 | 43 | 62 | 7 |
|---|---|---|---|---|---|---|---|---|---|---|---|---|---|---|---|---|---|---|
| 46 | 59 | 2 | 23 | | 17 | 12 | 61 | 40 | | 4 | 53 | 48 | 25 | | 63 | 6 | 19 | 42 |
| 31 | 38 | 51 | 10 | | 36 | 21 | 16 | 57 | | 49 | 44 | 29 | 8 | | 14 | 27 | 34 | 55 |
| 50 | 11 | 30 | 39 | | 13 | 60 | 33 | 24 | | 32 | 5 | 52 | 41 | | 35 | 54 | 15 | 26 |

D1= 70  D2= 198  D3= 126  D4= 126

### 60. Awani Kumar 2017

| 3 | 22 | 47 | 58 | | 64 | 5 | 52 | 9 | | 13 | 60 | 1 | 56 | | 50 | 43 | 30 | 7 |
|---|---|---|---|---|---|---|---|---|---|---|---|---|---|---|---|---|---|---|
| 46 | 59 | 2 | 23 | | 17 | 44 | 29 | 40 | | 36 | 21 | 48 | 25 | | 31 | 6 | 51 | 42 |
| 63 | 38 | 19 | 10 | | 4 | 53 | 16 | 57 | | 49 | 12 | 61 | 8 | | 14 | 27 | 34 | 55 |
| 18 | 11 | 62 | 39 | | 45 | 28 | 33 | 24 | | 32 | 37 | 20 | 41 | | 35 | 54 | 15 | 26 |

D1= 134  D2= 134  D3= 126  D4= 126

### 61. Awani Kumar 2011

| 3 | 22 | 63 | 42 | | 62 | 43 | 2 | 23 | | 55 | 34 | 11 | 30 | | 10 | 31 | 54 | 35 |
|---|---|---|---|---|---|---|---|---|---|---|---|---|---|---|---|---|---|---|
| 56 | 41 | 4 | 29 | | 1 | 32 | 53 | 44 | | 12 | 21 | 64 | 33 | | 61 | 36 | 9 | 24 |
| 13 | 28 | 49 | 40 | | 52 | 37 | 16 | 25 | | 57 | 48 | 5 | 20 | | 8 | 17 | 60 | 45 |
| 58 | 39 | 14 | 19 | | 15 | 18 | 59 | 38 | | 6 | 27 | 50 | 47 | | 51 | 46 | 7 | 26 |

D1= 66  D2= 194  D3= 194  D4= 66

### 62. Awani Kumar 2009

| 3 | 22 | 63 | 42 | | 64 | 41 | 4 | 21 | | 49 | 12 | 29 | 40 | | 14 | 55 | 34 | 27 |
|---|---|---|---|---|---|---|---|---|---|---|---|---|---|---|---|---|---|---|
| 62 | 43 | 2 | 23 | | 1 | 24 | 61 | 44 | | 32 | 37 | 52 | 9 | | 35 | 26 | 15 | 54 |
| 19 | 58 | 47 | 6 | | 48 | 5 | 20 | 57 | | 13 | 56 | 33 | 28 | | 50 | 11 | 30 | 39 |
| 46 | 7 | 18 | 59 | | 17 | 60 | 45 | 8 | | 36 | 25 | 16 | 53 | | 31 | 38 | 51 | 10 |

D1= 70  D2= 190  D3= 130  D4= 130

### 63. Awani Kumar 2017 — D1= 66  D2= 194  D3= 194  D4= 66

| 3 | 26 | 47 | 54 |   | 46 | 55 | 2 | 27 |   | 63 | 38 | 19 | 10 |   | 18 | 11 | 62 | 39 |
|---|---|---|---|---|---|---|---|---|---|---|---|---|---|---|---|---|---|---|
| 32 | 37 | 20 | 41 |   | 49 | 12 | 61 | 8 |   | 36 | 25 | 48 | 21 |   | 13 | 56 | 1 | 60 |
| 45 | 24 | 33 | 28 |   | 4 | 57 | 16 | 53 |   | 17 | 44 | 29 | 40 |   | 64 | 5 | 52 | 9 |
| 50 | 43 | 30 | 7 |   | 31 | 6 | 51 | 42 |   | 14 | 23 | 34 | 59 |   | 35 | 58 | 15 | 22 |

### 64. Awani Kumar 2007 — D1= 66  D2= 194  D3= 194  D4= 66

| 3 | 26 | 47 | 54 |   | 64 | 37 | 20 | 9 |   | 45 | 56 | 1 | 28 |   | 18 | 11 | 62 | 39 |
|---|---|---|---|---|---|---|---|---|---|---|---|---|---|---|---|---|---|---|
| 46 | 55 | 2 | 27 |   | 17 | 12 | 61 | 40 |   | 4 | 25 | 48 | 53 |   | 63 | 38 | 19 | 10 |
| 31 | 6 | 51 | 42 |   | 36 | 57 | 16 | 21 |   | 49 | 44 | 29 | 8 |   | 14 | 23 | 34 | 59 |
| 50 | 43 | 30 | 7 |   | 13 | 24 | 33 | 60 |   | 32 | 5 | 52 | 41 |   | 35 | 58 | 15 | 22 |

### 65. Awani Kumar 2017 — D1= 130  D2= 130  D3= 130  D4= 130

| 3 | 26 | 47 | 54 |   | 64 | 5 | 52 | 9 |   | 13 | 56 | 1 | 60 |   | 50 | 43 | 30 | 7 |
|---|---|---|---|---|---|---|---|---|---|---|---|---|---|---|---|---|---|---|
| 46 | 55 | 2 | 27 |   | 17 | 44 | 29 | 40 |   | 36 | 25 | 48 | 21 |   | 31 | 6 | 51 | 42 |
| 63 | 38 | 19 | 10 |   | 4 | 57 | 16 | 53 |   | 49 | 12 | 61 | 8 |   | 14 | 23 | 34 | 59 |
| 18 | 11 | 62 | 39 |   | 45 | 24 | 33 | 28 |   | 32 | 37 | 20 | 41 |   | 35 | 58 | 15 | 22 |

### 66. Francis Gaspalou 2009 — D1= 66  D2= 178  D3= 210  D4= 66

| 3 | 26 | 63 | 38 |   | 64 | 37 | 20 | 9 |   | 45 | 56 | 1 | 28 |   | 18 | 11 | 46 | 55 |
|---|---|---|---|---|---|---|---|---|---|---|---|---|---|---|---|---|---|---|
| 62 | 39 | 2 | 27 |   | 17 | 12 | 61 | 40 |   | 4 | 25 | 48 | 53 |   | 47 | 54 | 19 | 10 |
| 15 | 6 | 51 | 58 |   | 36 | 57 | 16 | 21 |   | 49 | 44 | 29 | 8 |   | 30 | 23 | 34 | 43 |
| 50 | 59 | 14 | 7 |   | 13 | 24 | 33 | 60 |   | 32 | 5 | 52 | 41 |   | 35 | 42 | 31 | 22 |

### 67. Francis Gaspalou 2009 — D1= 50  D2= 178  D3= 210  D4= 82

| 3 | 26 | 63 | 38 |   | 64 | 37 | 20 | 9 |   | 45 | 56 | 1 | 28 |   | 18 | 11 | 46 | 55 |
|---|---|---|---|---|---|---|---|---|---|---|---|---|---|---|---|---|---|---|
| 62 | 39 | 2 | 27 |   | 17 | 12 | 61 | 40 |   | 4 | 25 | 48 | 53 |   | 47 | 54 | 19 | 10 |
| 15 | 22 | 51 | 42 |   | 36 | 57 | 16 | 21 |   | 49 | 44 | 29 | 8 |   | 30 | 7 | 34 | 59 |
| 50 | 43 | 14 | 23 |   | 13 | 24 | 33 | 60 |   | 32 | 5 | 52 | 41 |   | 35 | 58 | 31 | 6 |

### 68. Francis Gaspalou 2009 — D1= 66  D2= 194  D3= 194  D4= 66

| 3 | 26 | 63 | 38 |   | 64 | 37 | 20 | 9 |   | 45 | 56 | 1 | 28 |   | 18 | 11 | 46 | 55 |
|---|---|---|---|---|---|---|---|---|---|---|---|---|---|---|---|---|---|---|
| 62 | 39 | 2 | 27 |   | 17 | 12 | 61 | 40 |   | 4 | 25 | 48 | 53 |   | 47 | 54 | 19 | 10 |
| 31 | 6 | 35 | 58 |   | 36 | 57 | 16 | 21 |   | 49 | 44 | 29 | 8 |   | 14 | 23 | 50 | 43 |
| 34 | 59 | 30 | 7 |   | 13 | 24 | 33 | 60 |   | 32 | 5 | 52 | 41 |   | 51 | 42 | 15 | 22 |

### 69. Francis Gaspalou 2009 — D1= 50  D2= 194  D3= 194  D4= 82

| 3 | 26 | 63 | 38 |   | 64 | 37 | 20 | 9 |   | 45 | 56 | 1 | 28 |   | 18 | 11 | 46 | 55 |
|---|---|---|---|---|---|---|---|---|---|---|---|---|---|---|---|---|---|---|
| 62 | 39 | 2 | 27 |   | 17 | 12 | 61 | 40 |   | 4 | 25 | 48 | 53 |   | 47 | 54 | 19 | 10 |
| 31 | 22 | 35 | 42 |   | 36 | 57 | 16 | 21 |   | 49 | 44 | 29 | 8 |   | 14 | 7 | 50 | 59 |
| 34 | 43 | 30 | 23 |   | 13 | 24 | 33 | 60 |   | 32 | 5 | 52 | 41 |   | 51 | 58 | 15 | 6 |

### 70. Francis Gaspalou 2009 — D1= 130  D2= 130  D3= 130  D4= 130

| 3 | 30 | 63 | 34 |   | 58 | 39 | 6 | 27 |   | 31 | 2 | 35 | 62 |   | 38 | 59 | 26 | 7 |
|---|---|---|---|---|---|---|---|---|---|---|---|---|---|---|---|---|---|---|
| 48 | 1 | 20 | 61 |   | 21 | 60 | 41 | 8 |   | 4 | 45 | 64 | 17 |   | 57 | 24 | 5 | 44 |
| 29 | 52 | 33 | 16 |   | 40 | 9 | 28 | 53 |   | 49 | 32 | 13 | 36 |   | 12 | 37 | 56 | 25 |
| 50 | 47 | 14 | 19 |   | 11 | 22 | 55 | 42 |   | 46 | 51 | 18 | 15 |   | 23 | 10 | 43 | 54 |

### 71. Awani Kumar 2017 — D1= 102  D2= 158  D3= 94  D4= 166

| 3 | 32 | 45 | 50 |   | 54 | 37 | 28 | 11 |   | 47 | 20 | 33 | 30 |   | 26 | 41 | 24 | 39 |
|---|---|---|---|---|---|---|---|---|---|---|---|---|---|---|---|---|---|---|
| 46 | 49 | 4 | 31 |   | 27 | 12 | 53 | 38 |   | 2 | 61 | 16 | 51 |   | 55 | 8 | 57 | 10 |
| 63 | 36 | 17 | 14 |   | 6 | 21 | 44 | 59 |   | 19 | 48 | 29 | 34 |   | 42 | 25 | 40 | 23 |
| 18 | 13 | 64 | 35 |   | 43 | 60 | 5 | 22 |   | 62 | 1 | 52 | 15 |   | 7 | 56 | 9 | 58 |

| 3 | 32 | 45 | 50 | | 58 | 37 | 24 | 11 | | 47 | 20 | 33 | 30 | | 22 | 41 | 28 | 39 |
|---|---|---|---|---|---|---|---|---|---|---|---|---|---|---|---|---|---|---|
| 46 | 49 | 4 | 31 | | 23 | 12 | 57 | 38 | | 2 | 61 | 16 | 51 | | 59 | 8 | 53 | 10 |
| 63 | 36 | 17 | 14 | | 6 | 25 | 44 | 55 | | 19 | 48 | 29 | 34 | | 42 | 21 | 40 | 27 |
| 18 | 13 | 64 | 35 | | 43 | 56 | 5 | 26 | | 62 | 1 | 52 | 15 | | 7 | 60 | 9 | 54 |
| 72. Awani Kumar 2012 | | | | | D1= | 98 | | D2= | 162 | | D3= | 98 | | D4= | 162 | | | |
| 3 | 42 | 63 | 22 | | 62 | 23 | 2 | 43 | | 47 | 6 | 19 | 58 | | 18 | 59 | 46 | 7 |
| 48 | 21 | 4 | 57 | | 1 | 60 | 45 | 24 | | 20 | 41 | 64 | 5 | | 61 | 8 | 17 | 44 |
| 29 | 40 | 49 | 12 | | 52 | 9 | 32 | 37 | | 33 | 28 | 13 | 56 | | 16 | 53 | 36 | 25 |
| 50 | 27 | 14 | 39 | | 15 | 38 | 51 | 26 | | 30 | 55 | 34 | 11 | | 35 | 10 | 31 | 54 |
| 73. Francis Gaspalou 2009 | | | | | D1= | 130 | | D2= | 130 | | D3= | 130 | | D4= | 130 | | | |
| 3 | 42 | 63 | 22 | | 62 | 23 | 2 | 43 | | 55 | 30 | 11 | 34 | | 10 | 35 | 54 | 31 |
| 56 | 29 | 4 | 41 | | 1 | 44 | 53 | 32 | | 12 | 33 | 64 | 21 | | 61 | 24 | 9 | 36 |
| 13 | 40 | 49 | 28 | | 52 | 25 | 16 | 37 | | 57 | 20 | 5 | 48 | | 8 | 45 | 60 | 17 |
| 58 | 19 | 14 | 39 | | 15 | 38 | 59 | 18 | | 6 | 47 | 50 | 27 | | 51 | 26 | 7 | 46 |
| 74. Awani Kumar 2012 | | | | | D1= | 98 | | D2= | 146 | | D3= | 178 | | D4= | 98 | | | |
| 3 | 42 | 63 | 22 | | 64 | 21 | 4 | 41 | | 29 | 40 | 49 | 12 | | 34 | 27 | 14 | 55 |
| 62 | 23 | 2 | 43 | | 1 | 44 | 61 | 24 | | 52 | 9 | 32 | 37 | | 15 | 54 | 35 | 26 |
| 47 | 6 | 19 | 58 | | 20 | 57 | 48 | 5 | | 33 | 28 | 13 | 56 | | 30 | 39 | 50 | 11 |
| 18 | 59 | 46 | 7 | | 45 | 8 | 17 | 60 | | 16 | 53 | 36 | 25 | | 51 | 10 | 31 | 38 |
| 75. Awani Kumar 2009 | | | | | D1= | 98 | | D2= | 162 | | D3= | 162 | | D4= | 98 | | | |
| 3 | 44 | 21 | 62 | | 54 | 29 | 4 | 43 | | 45 | 6 | 59 | 20 | | 28 | 51 | 46 | 5 |
| 48 | 7 | 58 | 17 | | 25 | 50 | 47 | 8 | | 2 | 41 | 24 | 63 | | 55 | 32 | 1 | 42 |
| 53 | 30 | 35 | 12 | | 36 | 11 | 22 | 61 | | 27 | 52 | 13 | 38 | | 14 | 37 | 60 | 19 |
| 26 | 49 | 16 | 39 | | 15 | 40 | 57 | 18 | | 56 | 31 | 34 | 9 | | 33 | 10 | 23 | 64 |
| 76. Guenter Stertenbrink 2003 | | | | | D1= | 130 | | D2= | 194 | | D3= | 66 | | D4= | 130 | | | |
| 3 | 46 | 31 | 50 | | 64 | 17 | 36 | 13 | | 45 | 4 | 49 | 32 | | 18 | 63 | 14 | 35 |
| 54 | 27 | 6 | 43 | | 37 | 12 | 21 | 60 | | 28 | 53 | 44 | 5 | | 11 | 38 | 59 | 22 |
| 47 | 2 | 51 | 30 | | 20 | 61 | 16 | 33 | | 1 | 48 | 29 | 52 | | 62 | 19 | 34 | 15 |
| 26 | 55 | 42 | 7 | | 9 | 40 | 57 | 24 | | 56 | 25 | 8 | 41 | | 39 | 10 | 23 | 58 |
| 77. Awani Kumar 2007 | | | | | D1= | 102 | | D2= | 158 | | D3= | 166 | | D4= | 94 | | | |
| 3 | 46 | 63 | 18 | | 58 | 23 | 6 | 43 | | 47 | 2 | 19 | 62 | | 22 | 59 | 42 | 7 |
| 48 | 1 | 20 | 61 | | 21 | 60 | 41 | 8 | | 4 | 45 | 64 | 17 | | 57 | 24 | 5 | 44 |
| 29 | 52 | 33 | 16 | | 40 | 9 | 28 | 53 | | 49 | 32 | 13 | 36 | | 12 | 37 | 56 | 25 |
| 50 | 31 | 14 | 35 | | 11 | 38 | 55 | 26 | | 30 | 51 | 34 | 15 | | 39 | 10 | 27 | 54 |
| 78. Awani Kumar 2007 | | | | | D1= | 130 | | D2= | 130 | | D3= | 130 | | D4= | 130 | | | |
| 3 | 46 | 63 | 18 | | 58 | 23 | 42 | 7 | | 47 | 2 | 19 | 62 | | 22 | 59 | 6 | 43 |
| 48 | 1 | 20 | 61 | | 21 | 60 | 5 | 44 | | 4 | 45 | 64 | 17 | | 57 | 24 | 41 | 8 |
| 29 | 52 | 33 | 16 | | 12 | 37 | 28 | 53 | | 49 | 32 | 13 | 36 | | 40 | 9 | 56 | 25 |
| 50 | 31 | 14 | 35 | | 39 | 10 | 55 | 26 | | 30 | 51 | 34 | 15 | | 11 | 38 | 27 | 54 |
| 79. Awani Kumar 2007 | | | | | D1= | 130 | | D2= | 66 | | D3= | 194 | | D4= | 130 | | | |
| 3 | 46 | 63 | 18 | | 58 | 23 | 6 | 43 | | 47 | 2 | 19 | 62 | | 22 | 59 | 42 | 7 |
| 48 | 17 | 4 | 61 | | 5 | 60 | 41 | 24 | | 20 | 45 | 64 | 1 | | 57 | 8 | 21 | 44 |
| 29 | 52 | 33 | 16 | | 40 | 9 | 28 | 53 | | 49 | 32 | 13 | 36 | | 12 | 37 | 56 | 25 |
| 50 | 15 | 30 | 35 | | 27 | 38 | 55 | 10 | | 14 | 51 | 34 | 31 | | 39 | 26 | 11 | 54 |
| 80. Awani Kumar 2007 | | | | | D1= | 130 | | D2= | 130 | | D3= | 130 | | D4= | 130 | | | |

| 3 | 46 | 63 | 18 |
|---|----|----|----|
| 54 | 27 | 6 | 43 |
| 47 | 2 | 19 | 62 |
| 26 | 55 | 42 | 7 |

| 64 | 17 | 4 | 45 |
|---|----|----|----|
| 5 | 44 | 53 | 28 |
| 52 | 29 | 16 | 33 |
| 9 | 40 | 57 | 24 |

| 13 | 36 | 49 | 32 |
|---|----|----|----|
| 60 | 21 | 12 | 37 |
| 1 | 48 | 61 | 20 |
| 56 | 25 | 8 | 41 |

| 50 | 31 | 14 | 35 |
|---|----|----|----|
| 11 | 38 | 59 | 22 |
| 30 | 51 | 34 | 15 |
| 39 | 10 | 23 | 58 |

81. Awani Kumar 2017    D1= 166    D2= 158    D3= 102    D4= 94

| 3 | 46 | 63 | 18 |
|---|----|----|----|
| 58 | 23 | 6 | 43 |
| 47 | 2 | 19 | 62 |
| 22 | 59 | 42 | 7 |

| 64 | 1 | 20 | 45 |
|---|----|----|----|
| 21 | 44 | 57 | 8 |
| 4 | 61 | 48 | 17 |
| 41 | 24 | 5 | 60 |

| 29 | 52 | 33 | 16 |
|---|----|----|----|
| 40 | 9 | 28 | 53 |
| 49 | 32 | 13 | 36 |
| 12 | 37 | 56 | 25 |

| 34 | 31 | 14 | 51 |
|---|----|----|----|
| 11 | 54 | 39 | 26 |
| 30 | 35 | 50 | 15 |
| 55 | 10 | 27 | 38 |

82. Awani Kumar 2007    D1= 98    D2= 162    D3= 162    D4= 98

| 3 | 46 | 63 | 18 |
|---|----|----|----|
| 58 | 23 | 6 | 43 |
| 47 | 2 | 19 | 62 |
| 22 | 59 | 42 | 7 |

| 64 | 17 | 4 | 45 |
|---|----|----|----|
| 5 | 44 | 57 | 24 |
| 20 | 61 | 48 | 1 |
| 41 | 8 | 21 | 60 |

| 29 | 52 | 33 | 16 |
|---|----|----|----|
| 40 | 9 | 28 | 53 |
| 49 | 32 | 13 | 36 |
| 12 | 37 | 56 | 25 |

| 34 | 15 | 30 | 51 |
|---|----|----|----|
| 27 | 54 | 39 | 10 |
| 14 | 35 | 50 | 31 |
| 55 | 26 | 11 | 38 |

83. Awani Kumar 2007    D1= 98    D2= 162    D3= 162    D4= 98

| 3 | 46 | 63 | 18 |
|---|----|----|----|
| 58 | 23 | 6 | 43 |
| 47 | 2 | 19 | 62 |
| 22 | 59 | 42 | 7 |

| 64 | 17 | 4 | 45 |
|---|----|----|----|
| 5 | 44 | 57 | 24 |
| 52 | 29 | 16 | 33 |
| 9 | 40 | 53 | 28 |

| 13 | 36 | 49 | 32 |
|---|----|----|----|
| 56 | 25 | 12 | 37 |
| 1 | 48 | 61 | 20 |
| 60 | 21 | 8 | 41 |

| 50 | 31 | 14 | 35 |
|---|----|----|----|
| 11 | 38 | 55 | 26 |
| 30 | 51 | 34 | 15 |
| 39 | 10 | 27 | 54 |

84. Awani Kumar 2012    D1= 162    D2= 162    D3= 98    D4= 98

| 3 | 48 | 29 | 50 |
|---|----|----|----|
| 58 | 21 | 40 | 11 |
| 63 | 4 | 49 | 14 |
| 6 | 57 | 12 | 55 |

| 62 | 1 | 52 | 15 |
|---|----|----|----|
| 7 | 60 | 9 | 54 |
| 2 | 45 | 32 | 51 |
| 59 | 24 | 37 | 10 |

| 47 | 20 | 33 | 30 |
|---|----|----|----|
| 22 | 41 | 28 | 39 |
| 19 | 64 | 13 | 34 |
| 42 | 5 | 56 | 27 |

| 18 | 61 | 16 | 35 |
|---|----|----|----|
| 43 | 8 | 53 | 26 |
| 46 | 17 | 36 | 31 |
| 23 | 44 | 25 | 38 |

85. Awani Kumar 2009    D1= 114    D2= 146    D3= 114    D4= 146

| 4 | 21 | 48 | 57 |
|---|----|----|----|
| 45 | 60 | 1 | 24 |
| 64 | 5 | 20 | 41 |
| 17 | 44 | 61 | 8 |

| 49 | 12 | 29 | 40 |
|---|----|----|----|
| 32 | 37 | 52 | 9 |
| 13 | 28 | 33 | 56 |
| 36 | 53 | 16 | 25 |

| 46 | 59 | 2 | 23 |
|---|----|----|----|
| 3 | 22 | 47 | 58 |
| 18 | 43 | 62 | 7 |
| 63 | 6 | 19 | 42 |

| 31 | 38 | 51 | 10 |
|---|----|----|----|
| 50 | 11 | 30 | 39 |
| 35 | 54 | 15 | 26 |
| 14 | 27 | 34 | 55 |

86. Awani Kumar 2006    D1= 158    D2= 166    D3= 102    D4= 94

| 4 | 21 | 62 | 43 |
|---|----|----|----|
| 59 | 46 | 5 | 20 |
| 22 | 3 | 44 | 61 |
| 45 | 60 | 19 | 6 |

| 63 | 42 | 1 | 24 |
|---|----|----|----|
| 8 | 17 | 58 | 47 |
| 41 | 64 | 23 | 2 |
| 18 | 7 | 48 | 57 |

| 54 | 35 | 12 | 29 |
|---|----|----|----|
| 13 | 28 | 51 | 38 |
| 36 | 53 | 30 | 11 |
| 27 | 14 | 37 | 52 |

| 9 | 32 | 55 | 34 |
|---|----|----|----|
| 50 | 39 | 16 | 25 |
| 31 | 10 | 33 | 56 |
| 40 | 49 | 26 | 15 |

87. Francis Gaspalou 2009    D1= 66    D2= 194    D3= 194    D4= 66

| 4 | 23 | 58 | 45 |
|---|----|----|----|
| 57 | 46 | 3 | 24 |
| 62 | 9 | 40 | 19 |
| 7 | 52 | 29 | 42 |

| 63 | 44 | 5 | 18 |
|---|----|----|----|
| 6 | 17 | 64 | 43 |
| 33 | 22 | 59 | 16 |
| 28 | 47 | 2 | 53 |

| 26 | 13 | 36 | 55 |
|---|----|----|----|
| 35 | 56 | 25 | 14 |
| 8 | 51 | 30 | 41 |
| 61 | 10 | 39 | 20 |

| 37 | 50 | 31 | 12 |
|---|----|----|----|
| 32 | 11 | 38 | 49 |
| 27 | 48 | 1 | 54 |
| 34 | 21 | 60 | 15 |

88. Awani Kumar 2017    D1= 66    D2= 194    D3= 66    D4= 194

| 4 | 23 | 58 | 45 |
|---|----|----|----|
| 57 | 46 | 3 | 24 |
| 62 | 41 | 8 | 19 |
| 7 | 20 | 61 | 42 |

| 63 | 44 | 5 | 18 |
|---|----|----|----|
| 6 | 17 | 64 | 43 |
| 1 | 22 | 59 | 48 |
| 60 | 47 | 2 | 21 |

| 10 | 29 | 52 | 39 |
|---|----|----|----|
| 51 | 40 | 9 | 30 |
| 56 | 35 | 14 | 25 |
| 13 | 26 | 55 | 36 |

| 53 | 34 | 15 | 28 |
|---|----|----|----|
| 16 | 27 | 54 | 33 |
| 11 | 32 | 49 | 38 |
| 50 | 37 | 12 | 31 |

89. Francis Gaspalou 2009    D1= 66    D2= 194    D3= 66    D4= 194

| 4 | 23 | 58 | 45 | | 63 | 44 | 5 | 18 | | 26 | 13 | 36 | 55 | | 37 | 50 | 31 | 12 |
|---|---|---|---|---|---|---|---|---|---|---|---|---|---|---|---|---|---|---|
| 57 | 46 | 3 | 24 | | 6 | 17 | 64 | 43 | | 35 | 56 | 25 | 14 | | 32 | 11 | 38 | 49 |
| 62 | 41 | 8 | 19 | | 1 | 22 | 59 | 48 | | 40 | 51 | 30 | 9 | | 27 | 16 | 33 | 54 |
| 7 | 20 | 61 | 42 | | 60 | 47 | 2 | 21 | | 29 | 10 | 39 | 52 | | 34 | 53 | 28 | 15 |

90. Francis Gaspalou 2009  D1= 66  D2= 194  D3= 66  D4= 194

| 4 | 25 | 46 | 55 | | 63 | 38 | 1 | 28 | | 26 | 35 | 56 | 13 | | 37 | 32 | 27 | 34 |
|---|---|---|---|---|---|---|---|---|---|---|---|---|---|---|---|---|---|---|
| 47 | 54 | 17 | 12 | | 20 | 9 | 62 | 39 | | 53 | 16 | 11 | 50 | | 10 | 51 | 40 | 29 |
| 58 | 3 | 24 | 45 | | 5 | 64 | 59 | 2 | | 36 | 57 | 14 | 23 | | 31 | 6 | 33 | 60 |
| 21 | 48 | 43 | 18 | | 42 | 19 | 8 | 61 | | 15 | 22 | 49 | 44 | | 52 | 41 | 30 | 7 |

91. Awani Kumar 2017  D1= 34  D2= 226  D3= 130  D4= 130

| 4 | 27 | 54 | 45 | | 63 | 40 | 9 | 18 | | 22 | 13 | 36 | 59 | | 41 | 50 | 31 | 8 |
|---|---|---|---|---|---|---|---|---|---|---|---|---|---|---|---|---|---|---|
| 53 | 46 | 3 | 28 | | 10 | 17 | 64 | 39 | | 35 | 60 | 21 | 14 | | 32 | 7 | 42 | 49 |
| 62 | 5 | 44 | 19 | | 33 | 26 | 55 | 16 | | 12 | 51 | 30 | 37 | | 23 | 48 | 1 | 58 |
| 11 | 52 | 29 | 38 | | 24 | 47 | 2 | 57 | | 61 | 6 | 43 | 20 | | 34 | 25 | 56 | 15 |

92. Awani Kumar 2017  D1= 66  D2= 194  D3= 66  D4= 194

| 4 | 29 | 58 | 39 | | 63 | 34 | 5 | 28 | | 10 | 55 | 52 | 13 | | 53 | 12 | 15 | 50 |
|---|---|---|---|---|---|---|---|---|---|---|---|---|---|---|---|---|---|---|
| 57 | 40 | 3 | 30 | | 6 | 27 | 64 | 33 | | 19 | 46 | 41 | 24 | | 48 | 17 | 22 | 43 |
| 62 | 35 | 8 | 25 | | 1 | 32 | 59 | 38 | | 56 | 9 | 14 | 51 | | 11 | 54 | 49 | 16 |
| 7 | 26 | 61 | 36 | | 60 | 37 | 2 | 31 | | 45 | 20 | 23 | 42 | | 18 | 47 | 44 | 21 |

93. Awani Kumar 2017  D1= 66  D2= 130  D3= 130  D4= 194

| 4 | 31 | 58 | 37 | | 63 | 36 | 5 | 26 | | 50 | 45 | 12 | 23 | | 13 | 18 | 55 | 44 |
|---|---|---|---|---|---|---|---|---|---|---|---|---|---|---|---|---|---|---|
| 49 | 46 | 3 | 32 | | 14 | 17 | 64 | 35 | | 11 | 24 | 57 | 38 | | 56 | 43 | 6 | 25 |
| 30 | 33 | 40 | 27 | | 1 | 62 | 59 | 8 | | 48 | 19 | 22 | 41 | | 51 | 16 | 9 | 54 |
| 47 | 20 | 29 | 34 | | 52 | 15 | 2 | 61 | | 21 | 42 | 39 | 28 | | 10 | 53 | 60 | 7 |

94. Awani Kumar 2017  D1= 50  D2= 130  D3= 210  D4= 130

| 4 | 37 | 62 | 27 | | 63 | 26 | 5 | 36 | | 54 | 3 | 12 | 61 | | 9 | 64 | 51 | 6 |
|---|---|---|---|---|---|---|---|---|---|---|---|---|---|---|---|---|---|---|
| 59 | 30 | 1 | 40 | | 8 | 33 | 58 | 31 | | 13 | 60 | 55 | 2 | | 50 | 7 | 16 | 57 |
| 38 | 19 | 28 | 45 | | 25 | 48 | 35 | 22 | | 20 | 53 | 46 | 11 | | 47 | 10 | 21 | 52 |
| 29 | 44 | 39 | 18 | | 34 | 23 | 32 | 41 | | 43 | 14 | 17 | 56 | | 24 | 49 | 42 | 15 |

95. Awani Kumar 2012  D1= 98  D2= 162  D3= 138  D4= 122

| 4 | 45 | 56 | 25 | | 57 | 24 | 5 | 44 | | 46 | 3 | 26 | 55 | | 23 | 58 | 43 | 6 |
|---|---|---|---|---|---|---|---|---|---|---|---|---|---|---|---|---|---|---|
| 53 | 28 | 1 | 48 | | 16 | 33 | 52 | 29 | | 27 | 54 | 47 | 2 | | 34 | 15 | 30 | 51 |
| 10 | 39 | 62 | 19 | | 21 | 60 | 41 | 8 | | 40 | 9 | 20 | 61 | | 59 | 22 | 7 | 42 |
| 63 | 18 | 11 | 38 | | 36 | 13 | 32 | 49 | | 17 | 64 | 37 | 12 | | 14 | 35 | 50 | 31 |

96. Awani Kumar 2012  D1= 88  D2= 100  D3= 176  D4= 156

| 4 | 45 | 64 | 17 | | 57 | 24 | 5 | 44 | | 48 | 1 | 20 | 61 | | 21 | 60 | 41 | 8 |
|---|---|---|---|---|---|---|---|---|---|---|---|---|---|---|---|---|---|---|
| 49 | 32 | 13 | 36 | | 12 | 37 | 56 | 25 | | 29 | 52 | 33 | 16 | | 40 | 9 | 28 | 53 |
| 46 | 3 | 18 | 63 | | 23 | 58 | 43 | 6 | | 2 | 47 | 62 | 19 | | 59 | 22 | 7 | 42 |
| 31 | 50 | 35 | 14 | | 38 | 11 | 26 | 55 | | 51 | 30 | 15 | 34 | | 10 | 39 | 54 | 27 |

97. Awani Kumar 2006  D1= 130  D2= 130  D3= 130  D4= 130

| 5 | 16 | 59 | 50 | | 58 | 51 | 8 | 13 | | 15 | 6 | 49 | 60 | | 52 | 57 | 14 | 7 |
|---|---|---|---|---|---|---|---|---|---|---|---|---|---|---|---|---|---|---|
| 18 | 27 | 48 | 37 | | 45 | 40 | 19 | 26 | | 28 | 17 | 38 | 47 | | 39 | 46 | 25 | 20 |
| 63 | 54 | 1 | 12 | | 4 | 9 | 62 | 55 | | 53 | 64 | 11 | 2 | | 10 | 3 | 56 | 61 |
| 44 | 33 | 22 | 31 | | 23 | 30 | 41 | 36 | | 34 | 43 | 32 | 21 | | 29 | 24 | 35 | 42 |

98. Awani Kumar 2009  D1= 98  D2= 162  D3= 98  D4= 162

| | | | | | | | | | | | | | | | | |
|---|---|---|---|---|---|---|---|---|---|---|---|---|---|---|---|---|
| 5 | 18 | 63 | 44 | | 58 | 45 | 4 | 23 | | 31 | 12 | 37 | 50 | | 36 | 55 | 26 | 13 |
| 32 | 11 | 38 | 49 | | 35 | 56 | 25 | 14 | | 6 | 17 | 64 | 43 | | 57 | 46 | 3 | 24 |
| 59 | 48 | 1 | 22 | | 8 | 19 | 62 | 41 | | 33 | 54 | 27 | 16 | | 30 | 9 | 40 | 51 |
| 34 | 53 | 28 | 15 | | 29 | 10 | 39 | 52 | | 60 | 47 | 2 | 21 | | 7 | 20 | 61 | 42 |

99. Guenter Stertenbrink 2003    D1= 130    D2= 130    D3= 130    D4= 130

| 5 | 26 | 63 | 36 | | 64 | 35 | 6 | 25 | | 27 | 8 | 33 | 62 | | 34 | 61 | 28 | 7 |
|---|---|---|---|---|---|---|---|---|---|---|---|---|---|---|---|---|---|---|
| 50 | 45 | 12 | 23 | | 11 | 24 | 49 | 46 | | 48 | 51 | 22 | 9 | | 21 | 10 | 47 | 52 |
| 31 | 4 | 37 | 58 | | 38 | 57 | 32 | 3 | | 1 | 30 | 59 | 40 | | 60 | 39 | 2 | 29 |
| 44 | 55 | 18 | 13 | | 17 | 14 | 43 | 56 | | 54 | 41 | 16 | 19 | | 15 | 20 | 53 | 42 |

100. Guenter Stertenbrink 2003    D1= 130    D2= 130    D3= 130    D4= 130

| 5 | 26 | 63 | 36 | | 64 | 35 | 6 | 25 | | 27 | 40 | 1 | 62 | | 34 | 29 | 60 | 7 |
|---|---|---|---|---|---|---|---|---|---|---|---|---|---|---|---|---|---|---|
| 50 | 45 | 12 | 23 | | 11 | 24 | 49 | 46 | | 16 | 51 | 22 | 41 | | 53 | 10 | 47 | 20 |
| 31 | 4 | 37 | 58 | | 38 | 57 | 32 | 3 | | 33 | 30 | 59 | 8 | | 28 | 39 | 2 | 61 |
| 44 | 55 | 18 | 13 | | 17 | 14 | 43 | 56 | | 54 | 9 | 48 | 19 | | 15 | 52 | 21 | 42 |

101. Awani Kumar 2017    D1= 130    D2= 130    D3= 130    D4= 130

| 5 | 26 | 63 | 36 | | 58 | 37 | 4 | 31 | | 55 | 44 | 13 | 18 | | 12 | 23 | 50 | 45 |
|---|---|---|---|---|---|---|---|---|---|---|---|---|---|---|---|---|---|---|
| 56 | 43 | 6 | 25 | | 11 | 24 | 57 | 38 | | 14 | 17 | 64 | 35 | | 49 | 46 | 3 | 32 |
| 59 | 8 | 1 | 62 | | 40 | 27 | 30 | 33 | | 9 | 54 | 51 | 16 | | 22 | 41 | 48 | 19 |
| 10 | 53 | 60 | 7 | | 21 | 42 | 39 | 28 | | 52 | 15 | 2 | 61 | | 47 | 20 | 29 | 34 |

102. Awani Kumar 2017    D1= 114    D2= 194    D3= 146    D4= 66

| 5 | 32 | 59 | 34 | | 58 | 35 | 24 | 13 | | 31 | 38 | 33 | 28 | | 36 | 25 | 14 | 55 |
|---|---|---|---|---|---|---|---|---|---|---|---|---|---|---|---|---|---|---|
| 42 | 51 | 8 | 29 | | 21 | 16 | 43 | 50 | | 52 | 9 | 30 | 39 | | 15 | 54 | 49 | 12 |
| 63 | 6 | 1 | 60 | | 4 | 57 | 46 | 23 | | 37 | 64 | 27 | 2 | | 26 | 3 | 56 | 45 |
| 20 | 41 | 62 | 7 | | 47 | 22 | 17 | 44 | | 10 | 19 | 40 | 61 | | 53 | 48 | 11 | 18 |

103. Awani Kumar 2017    D1= 66    D2= 194    D3= 162    D4= 98

| 5 | 42 | 19 | 64 | | 52 | 31 | 6 | 41 | | 43 | 8 | 61 | 18 | | 30 | 49 | 44 | 7 |
|---|---|---|---|---|---|---|---|---|---|---|---|---|---|---|---|---|---|---|
| 46 | 1 | 60 | 23 | | 27 | 56 | 45 | 2 | | 4 | 47 | 22 | 57 | | 53 | 26 | 3 | 48 |
| 51 | 32 | 37 | 10 | | 38 | 9 | 20 | 63 | | 29 | 50 | 11 | 40 | | 12 | 39 | 62 | 17 |
| 28 | 55 | 14 | 33 | | 13 | 34 | 59 | 24 | | 54 | 25 | 36 | 15 | | 35 | 16 | 21 | 58 |

104. Guenter Stertenbrink 2003    D1= 130    D2= 194    D3= 66    D4= 130

| 6 | 19 | 58 | 47 | | 59 | 46 | 7 | 18 | | 50 | 39 | 14 | 27 | | 15 | 26 | 51 | 38 |
|---|---|---|---|---|---|---|---|---|---|---|---|---|---|---|---|---|---|---|
| 57 | 48 | 5 | 20 | | 8 | 17 | 60 | 45 | | 13 | 28 | 49 | 40 | | 52 | 37 | 16 | 25 |
| 4 | 29 | 56 | 41 | | 53 | 44 | 1 | 32 | | 64 | 33 | 12 | 21 | | 9 | 24 | 61 | 36 |
| 63 | 34 | 11 | 22 | | 10 | 23 | 62 | 35 | | 3 | 30 | 55 | 42 | | 54 | 43 | 2 | 31 |

105. Awani Kumar 2011    D1= 66    D2= 194    D3= 194    D4= 66

| 6 | 25 | 44 | 55 | | 51 | 48 | 29 | 2 | | 42 | 19 | 8 | 61 | | 31 | 38 | 49 | 12 |
|---|---|---|---|---|---|---|---|---|---|---|---|---|---|---|---|---|---|---|
| 43 | 56 | 5 | 26 | | 30 | 1 | 52 | 47 | | 7 | 62 | 41 | 20 | | 50 | 11 | 32 | 37 |
| 58 | 45 | 24 | 3 | | 15 | 28 | 33 | 54 | | 22 | 39 | 60 | 9 | | 35 | 18 | 13 | 64 |
| 23 | 4 | 57 | 46 | | 34 | 53 | 16 | 27 | | 59 | 10 | 21 | 40 | | 14 | 63 | 36 | 17 |

106. Awani Kumar 2011    D1= 84    D2= 160    D3= 104    D4= 172

| 6 | 35 | 64 | 25 | | 61 | 28 | 3 | 38 | | 36 | 21 | 26 | 47 | | 27 | 46 | 37 | 20 |
|---|---|---|---|---|---|---|---|---|---|---|---|---|---|---|---|---|---|---|
| 57 | 32 | 7 | 34 | | 2 | 39 | 60 | 29 | | 31 | 42 | 33 | 24 | | 40 | 17 | 30 | 43 |
| 52 | 5 | 10 | 63 | | 11 | 62 | 53 | 4 | | 22 | 51 | 48 | 9 | | 45 | 12 | 19 | 54 |
| 15 | 58 | 49 | 8 | | 56 | 1 | 14 | 59 | | 41 | 16 | 23 | 50 | | 18 | 55 | 44 | 13 |

107. Awani Kumar 2012    D1= 106    D2= 154    D3= 130    D4= 130

| 6 | 41 | 24 | 59 |
|---|---|---|---|
| 43 | 8 | 57 | 22 |
| 54 | 25 | 40 | 11 |
| 27 | 56 | 9 | 38 |

| 51 | 32 | 13 | 34 |
|---|---|---|---|
| 30 | 49 | 36 | 15 |
| 47 | 4 | 17 | 62 |
| 2 | 45 | 64 | 19 |

| 42 | 5 | 60 | 23 |
|---|---|---|---|
| 7 | 44 | 21 | 58 |
| 26 | 53 | 12 | 39 |
| 55 | 28 | 37 | 10 |

| 31 | 52 | 33 | 14 |
|---|---|---|---|
| 50 | 29 | 16 | 35 |
| 3 | 48 | 61 | 18 |
| 46 | 1 | 20 | 63 |

108. Awani Kumar 2007  D1= 130  D2= 194  D3= 66  D4= 130

| 6 | 47 | 58 | 19 |
|---|---|---|---|
| 57 | 20 | 5 | 48 |
| 4 | 41 | 56 | 29 |
| 63 | 22 | 11 | 34 |

| 59 | 18 | 7 | 46 |
|---|---|---|---|
| 8 | 45 | 60 | 17 |
| 53 | 32 | 1 | 44 |
| 10 | 35 | 62 | 23 |

| 50 | 27 | 14 | 39 |
|---|---|---|---|
| 13 | 40 | 49 | 28 |
| 64 | 21 | 12 | 33 |
| 3 | 42 | 55 | 30 |

| 15 | 38 | 51 | 26 |
|---|---|---|---|
| 52 | 25 | 16 | 37 |
| 9 | 36 | 61 | 24 |
| 54 | 31 | 2 | 43 |

109. Awani Kumar 2012  D1= 106  D2= 154  D3= 170  D4= 90

| 6 | 47 | 58 | 19 |
|---|---|---|---|
| 57 | 20 | 5 | 48 |
| 12 | 49 | 40 | 29 |
| 55 | 14 | 27 | 34 |

| 59 | 18 | 7 | 46 |
|---|---|---|---|
| 8 | 45 | 60 | 17 |
| 37 | 32 | 9 | 52 |
| 26 | 35 | 54 | 15 |

| 42 | 63 | 22 | 3 |
|---|---|---|---|
| 21 | 4 | 41 | 64 |
| 56 | 13 | 28 | 33 |
| 11 | 50 | 39 | 30 |

| 23 | 2 | 43 | 62 |
|---|---|---|---|
| 44 | 61 | 24 | 1 |
| 25 | 36 | 53 | 16 |
| 38 | 31 | 10 | 51 |

110. Awani Kumar 2009  D1= 130  D2= 130  D3= 190  D4= 70

| 6 | 51 | 60 | 13 |
|---|---|---|---|
| 57 | 16 | 7 | 50 |
| 52 | 5 | 14 | 59 |
| 15 | 58 | 49 | 8 |

| 61 | 12 | 3 | 54 |
|---|---|---|---|
| 2 | 55 | 64 | 9 |
| 11 | 62 | 53 | 4 |
| 56 | 1 | 10 | 63 |

| 20 | 37 | 46 | 27 |
|---|---|---|---|
| 47 | 26 | 17 | 40 |
| 38 | 19 | 28 | 45 |
| 25 | 48 | 39 | 18 |

| 43 | 30 | 21 | 36 |
|---|---|---|---|
| 24 | 33 | 42 | 31 |
| 29 | 44 | 35 | 22 |
| 34 | 23 | 32 | 41 |

111. Francis Gaspalou 2009  D1= 130  D2= 130  D3= 130  D4= 130

| 7 | 18 | 43 | 62 |
|---|---|---|---|
| 24 | 45 | 28 | 33 |
| 41 | 32 | 37 | 20 |
| 58 | 35 | 22 | 15 |

| 42 | 63 | 6 | 19 |
|---|---|---|---|
| 57 | 4 | 53 | 16 |
| 8 | 49 | 12 | 61 |
| 23 | 14 | 59 | 34 |

| 55 | 46 | 27 | 2 |
|---|---|---|---|
| 40 | 17 | 44 | 29 |
| 25 | 36 | 21 | 48 |
| 10 | 31 | 38 | 51 |

| 26 | 3 | 54 | 47 |
|---|---|---|---|
| 9 | 64 | 5 | 52 |
| 56 | 13 | 60 | 1 |
| 39 | 50 | 11 | 30 |

112. Awani Kumar 2017  D1= 62  D2= 190  D3= 198  D4= 70

| 7 | 18 | 43 | 62 |
|---|---|---|---|
| 42 | 63 | 6 | 19 |
| 55 | 14 | 27 | 34 |
| 26 | 35 | 54 | 15 |

| 56 | 45 | 28 | 1 |
|---|---|---|---|
| 25 | 4 | 53 | 48 |
| 8 | 49 | 44 | 29 |
| 41 | 32 | 5 | 52 |

| 9 | 64 | 37 | 20 |
|---|---|---|---|
| 40 | 17 | 12 | 61 |
| 57 | 36 | 21 | 16 |
| 24 | 13 | 60 | 33 |

| 58 | 3 | 22 | 47 |
|---|---|---|---|
| 23 | 46 | 59 | 2 |
| 10 | 31 | 38 | 51 |
| 39 | 50 | 11 | 30 |

113. Francis Gaspalou 2009  D1= 62  D2= 190  D3= 134  D4= 134

| 7 | 18 | 43 | 62 |
|---|---|---|---|
| 42 | 63 | 6 | 19 |
| 55 | 46 | 27 | 2 |
| 26 | 3 | 54 | 47 |

| 56 | 13 | 60 | 1 |
|---|---|---|---|
| 25 | 36 | 21 | 48 |
| 8 | 49 | 12 | 61 |
| 41 | 32 | 37 | 20 |

| 9 | 64 | 5 | 52 |
|---|---|---|---|
| 40 | 17 | 44 | 29 |
| 57 | 4 | 53 | 16 |
| 24 | 45 | 28 | 33 |

| 58 | 35 | 22 | 15 |
|---|---|---|---|
| 23 | 14 | 59 | 34 |
| 10 | 31 | 38 | 51 |
| 39 | 50 | 11 | 30 |

114. Awani Kumar 2017  D1= 126  D2= 126  D3= 134  D4= 134

| 7 | 24 | 41 | 58 |
|---|---|---|---|
| 42 | 57 | 8 | 23 |
| 55 | 40 | 25 | 10 |
| 26 | 9 | 56 | 39 |

| 50 | 13 | 32 | 35 |
|---|---|---|---|
| 31 | 36 | 49 | 14 |
| 46 | 17 | 4 | 63 |
| 3 | 64 | 45 | 18 |

| 43 | 60 | 5 | 22 |
|---|---|---|---|
| 6 | 21 | 44 | 59 |
| 27 | 12 | 53 | 38 |
| 54 | 37 | 28 | 11 |

| 30 | 33 | 52 | 15 |
|---|---|---|---|
| 51 | 16 | 29 | 34 |
| 2 | 61 | 48 | 19 |
| 47 | 20 | 1 | 62 |

115. Awani Kumar 2007  D1= 158  D2= 166  D3= 102  D4= 94

| 7 | 24 | 41 | 58 |
|---|---|---|---|
| 42 | 57 | 8 | 23 |
| 55 | 40 | 25 | 10 |
| 26 | 9 | 56 | 39 |

| 50 | 45 | 32 | 3 |
|---|---|---|---|
| 31 | 4 | 49 | 46 |
| 14 | 17 | 36 | 63 |
| 35 | 64 | 13 | 18 |

| 43 | 28 | 37 | 22 |
|---|---|---|---|
| 6 | 53 | 12 | 59 |
| 27 | 44 | 21 | 38 |
| 54 | 5 | 60 | 11 |

| 30 | 33 | 20 | 47 |
|---|---|---|---|
| 51 | 16 | 61 | 2 |
| 34 | 29 | 48 | 19 |
| 15 | 52 | 1 | 62 |

116. Awani Kumar 2017  D1= 94  D2= 166  D3= 102  D4= 158

| | | | | | | | | | | | | | | | |
|---|---|---|---|---|---|---|---|---|---|---|---|---|---|---|---|
| 7 | 26 | 51 | 46 | 52 | 37 | 16 | 25 | 9 | 24 | 61 | 36 | 62 | 43 | 2 | 23 |
| 50 | 47 | 6 | 27 | 13 | 28 | 49 | 40 | 64 | 33 | 12 | 21 | 3 | 22 | 63 | 42 |
| 59 | 38 | 15 | 18 | 8 | 17 | 60 | 45 | 53 | 44 | 1 | 32 | 10 | 31 | 54 | 35 |
| 14 | 19 | 58 | 39 | 57 | 48 | 5 | 20 | 4 | 29 | 56 | 41 | 55 | 34 | 11 | 30 |

117. Awani Kumar 2011    D1= 66    D2= 194    D3= 66    D4= 194

| | | | | | | | | | | | | | | | |
|---|---|---|---|---|---|---|---|---|---|---|---|---|---|---|---|
| 7 | 42 | 55 | 26 | 56 | 25 | 8 | 41 | 9 | 40 | 57 | 24 | 58 | 23 | 10 | 39 |
| 50 | 31 | 14 | 35 | 13 | 36 | 49 | 32 | 64 | 17 | 4 | 45 | 3 | 46 | 63 | 18 |
| 43 | 6 | 27 | 54 | 60 | 21 | 12 | 37 | 5 | 44 | 53 | 28 | 22 | 59 | 38 | 11 |
| 30 | 51 | 34 | 15 | 1 | 48 | 61 | 20 | 52 | 29 | 16 | 33 | 47 | 2 | 19 | 62 |

118. Awani Kumar 2017    D1= 158    D2= 166    D3= 94    D4= 102

| | | | | | | | | | | | | | | | |
|---|---|---|---|---|---|---|---|---|---|---|---|---|---|---|---|
| 7 | 42 | 55 | 26 | 58 | 23 | 10 | 39 | 41 | 8 | 25 | 56 | 24 | 57 | 40 | 9 |
| 50 | 31 | 46 | 3 | 35 | 14 | 63 | 18 | 32 | 49 | 4 | 45 | 13 | 36 | 17 | 64 |
| 43 | 6 | 27 | 54 | 22 | 59 | 38 | 11 | 5 | 44 | 53 | 28 | 60 | 21 | 12 | 37 |
| 30 | 51 | 2 | 47 | 15 | 34 | 19 | 62 | 52 | 29 | 48 | 1 | 33 | 16 | 61 | 20 |

119. Awani Kumar 2006    D1= 94    D2= 166    D3= 102    D4= 158

| | | | | | | | | | | | | | | | |
|---|---|---|---|---|---|---|---|---|---|---|---|---|---|---|---|
| 7 | 46 | 51 | 26 | 52 | 25 | 16 | 37 | 9 | 36 | 61 | 24 | 62 | 23 | 2 | 43 |
| 50 | 27 | 6 | 47 | 13 | 40 | 49 | 28 | 64 | 21 | 12 | 33 | 3 | 42 | 63 | 22 |
| 59 | 18 | 15 | 38 | 8 | 45 | 60 | 17 | 53 | 32 | 1 | 44 | 10 | 35 | 54 | 31 |
| 14 | 39 | 58 | 19 | 57 | 20 | 5 | 48 | 4 | 41 | 56 | 29 | 55 | 30 | 11 | 34 |

120. Awani Kumar 2012    D1= 82    D2= 162    D3= 114    D4= 162

| | | | | | | | | | | | | | | | |
|---|---|---|---|---|---|---|---|---|---|---|---|---|---|---|---|
| 7 | 46 | 51 | 26 | 60 | 17 | 8 | 45 | 9 | 36 | 61 | 24 | 54 | 31 | 10 | 35 |
| 50 | 27 | 6 | 47 | 5 | 48 | 57 | 20 | 64 | 21 | 12 | 33 | 11 | 34 | 55 | 30 |
| 59 | 18 | 15 | 38 | 16 | 37 | 52 | 25 | 53 | 32 | 1 | 44 | 2 | 43 | 62 | 23 |
| 14 | 39 | 58 | 19 | 49 | 28 | 13 | 40 | 4 | 41 | 56 | 29 | 63 | 22 | 3 | 42 |

121. Awani Kumar 2012    D1= 98    D2= 178    D3= 98    D4= 146

| | | | | | | | | | | | | | | | |
|---|---|---|---|---|---|---|---|---|---|---|---|---|---|---|---|
| 7 | 46 | 57 | 20 | 60 | 17 | 6 | 47 | 13 | 56 | 51 | 10 | 50 | 11 | 16 | 53 |
| 58 | 19 | 8 | 45 | 5 | 48 | 59 | 18 | 36 | 25 | 30 | 39 | 31 | 38 | 33 | 28 |
| 61 | 24 | 3 | 42 | 2 | 43 | 64 | 21 | 55 | 14 | 9 | 52 | 12 | 49 | 54 | 15 |
| 4 | 41 | 62 | 23 | 63 | 22 | 1 | 44 | 26 | 35 | 40 | 29 | 37 | 32 | 27 | 34 |

122. Awani Kumar 2012    D1= 98    D2= 130    D3= 130    D4= 162

| | | | | | | | | | | | | | | | |
|---|---|---|---|---|---|---|---|---|---|---|---|---|---|---|---|
| 7 | 46 | 59 | 18 | 60 | 17 | 8 | 45 | 9 | 36 | 61 | 24 | 54 | 31 | 2 | 43 |
| 58 | 19 | 6 | 47 | 5 | 48 | 57 | 20 | 64 | 21 | 12 | 33 | 3 | 42 | 55 | 30 |
| 51 | 26 | 15 | 38 | 16 | 37 | 52 | 25 | 53 | 32 | 1 | 44 | 10 | 35 | 62 | 23 |
| 14 | 39 | 50 | 27 | 49 | 28 | 13 | 40 | 4 | 41 | 56 | 29 | 63 | 22 | 11 | 34 |

123. Awani Kumar 2012    D1= 90    D2= 170    D3= 106    D4= 154

| | | | | | | | | | | | | | | | |
|---|---|---|---|---|---|---|---|---|---|---|---|---|---|---|---|
| 8 | 17 | 44 | 61 | 41 | 64 | 5 | 20 | 24 | 45 | 60 | 1 | 57 | 4 | 21 | 48 |
| 25 | 36 | 53 | 16 | 56 | 13 | 28 | 33 | 9 | 32 | 37 | 52 | 40 | 49 | 12 | 29 |
| 42 | 63 | 6 | 19 | 7 | 18 | 43 | 62 | 58 | 3 | 22 | 47 | 23 | 46 | 59 | 2 |
| 55 | 14 | 27 | 34 | 26 | 35 | 54 | 15 | 39 | 50 | 11 | 30 | 10 | 31 | 38 | 51 |

124. Awani Kumar 2006    D1= 94    D2= 102    D3= 158    D4= 166

| | | | | | | | | | | | | | | | |
|---|---|---|---|---|---|---|---|---|---|---|---|---|---|---|---|
| 9 | 22 | 63 | 36 | 64 | 35 | 10 | 21 | 23 | 44 | 1 | 62 | 34 | 29 | 56 | 11 |
| 50 | 45 | 8 | 27 | 7 | 28 | 49 | 46 | 16 | 51 | 26 | 37 | 57 | 6 | 47 | 20 |
| 31 | 4 | 41 | 54 | 42 | 53 | 32 | 3 | 33 | 30 | 55 | 12 | 24 | 43 | 2 | 61 |
| 40 | 59 | 18 | 13 | 17 | 14 | 39 | 60 | 58 | 5 | 48 | 19 | 15 | 52 | 25 | 38 |

125. Awani Kumar 2017    D1= 130    D2= 130    D3= 130    D4= 130

| | | | | | | | | | | | | | | | |
|---|---|---|---|---|---|---|---|---|---|---|---|---|---|---|---|
| 10 | 3 | 54 | 63 | 55 | 62 | 11 | 2 | 6 | 15 | 58 | 51 | 59 | 50 | 7 | 14 |
| 37 | 64 | 9 | 20 | 12 | 17 | 40 | 61 | 57 | 36 | 21 | 16 | 24 | 13 | 60 | 33 |
| 56 | 45 | 28 | 1 | 25 | 4 | 53 | 48 | 44 | 49 | 8 | 29 | 5 | 32 | 41 | 52 |
| 27 | 18 | 39 | 46 | 38 | 47 | 26 | 19 | 23 | 30 | 43 | 34 | 42 | 35 | 22 | 31 |

126. Awani Kumar 2009 — D1= 66  D2= 194  D3= 66  D4= 194

| | | | | | | | | | | | | | | | |
|---|---|---|---|---|---|---|---|---|---|---|---|---|---|---|---|
| 10 | 3 | 54 | 63 | 55 | 62 | 11 | 2 | 6 | 31 | 58 | 35 | 59 | 34 | 7 | 30 |
| 37 | 64 | 9 | 20 | 12 | 17 | 40 | 61 | 57 | 36 | 21 | 16 | 24 | 13 | 60 | 33 |
| 56 | 45 | 28 | 1 | 25 | 4 | 53 | 48 | 44 | 49 | 8 | 29 | 5 | 32 | 41 | 52 |
| 27 | 18 | 39 | 46 | 38 | 47 | 26 | 19 | 23 | 14 | 43 | 50 | 42 | 51 | 22 | 15 |

127. Francis Gaspalou 2009 — D1= 50  D2= 194  D3= 82  D4= 194

| | | | | | | | | | | | | | | | |
|---|---|---|---|---|---|---|---|---|---|---|---|---|---|---|---|
| 10 | 3 | 54 | 63 | 55 | 62 | 11 | 2 | 22 | 15 | 42 | 51 | 43 | 50 | 23 | 14 |
| 37 | 64 | 9 | 20 | 12 | 17 | 40 | 61 | 57 | 36 | 21 | 16 | 24 | 13 | 60 | 33 |
| 56 | 45 | 28 | 1 | 25 | 4 | 53 | 48 | 44 | 49 | 8 | 29 | 5 | 32 | 41 | 52 |
| 27 | 18 | 39 | 46 | 38 | 47 | 26 | 19 | 7 | 30 | 59 | 34 | 58 | 35 | 6 | 31 |

128. Francis Gaspalou 2009 — D1= 66  D2= 210  D3= 66  D4= 178

| | | | | | | | | | | | | | | | |
|---|---|---|---|---|---|---|---|---|---|---|---|---|---|---|---|
| 10 | 3 | 54 | 63 | 55 | 62 | 11 | 2 | 22 | 31 | 42 | 35 | 43 | 34 | 23 | 30 |
| 37 | 64 | 9 | 20 | 12 | 17 | 40 | 61 | 57 | 36 | 21 | 16 | 24 | 13 | 60 | 33 |
| 56 | 45 | 28 | 1 | 25 | 4 | 53 | 48 | 44 | 49 | 8 | 29 | 5 | 32 | 41 | 52 |
| 27 | 18 | 39 | 46 | 38 | 47 | 26 | 19 | 7 | 14 | 59 | 50 | 58 | 51 | 6 | 15 |

129. Francis Gaspalou 2009 — D1= 50  D2= 210  D3= 82  D4= 178

| | | | | | | | | | | | | | | | |
|---|---|---|---|---|---|---|---|---|---|---|---|---|---|---|---|
| 10 | 3 | 56 | 61 | 53 | 64 | 11 | 2 | 4 | 9 | 62 | 55 | 63 | 54 | 1 | 12 |
| 29 | 24 | 35 | 42 | 34 | 43 | 32 | 21 | 23 | 30 | 41 | 36 | 44 | 33 | 22 | 31 |
| 52 | 57 | 14 | 7 | 15 | 6 | 49 | 60 | 58 | 51 | 8 | 13 | 5 | 16 | 59 | 50 |
| 39 | 46 | 25 | 20 | 28 | 17 | 38 | 47 | 45 | 40 | 19 | 26 | 18 | 27 | 48 | 37 |

130. Awani Kumar 2009 — D1= 98  D2= 162  D3= 98  D4= 162

| | | | | | | | | | | | | | | | |
|---|---|---|---|---|---|---|---|---|---|---|---|---|---|---|---|
| 10 | 3 | 56 | 61 | 53 | 64 | 11 | 2 | 4 | 9 | 62 | 55 | 63 | 54 | 1 | 12 |
| 39 | 46 | 25 | 20 | 28 | 17 | 38 | 47 | 45 | 40 | 19 | 26 | 18 | 27 | 48 | 37 |
| 52 | 57 | 14 | 7 | 15 | 6 | 49 | 60 | 58 | 51 | 8 | 13 | 5 | 16 | 59 | 50 |
| 29 | 24 | 35 | 42 | 34 | 43 | 32 | 21 | 23 | 30 | 41 | 36 | 44 | 33 | 22 | 31 |

131. Awani Kumar 2009 — D1= 66  D2= 194  D3= 66  D4= 194

| | | | | | | | | | | | | | | | |
|---|---|---|---|---|---|---|---|---|---|---|---|---|---|---|---|
| 10 | 19 | 38 | 63 | 39 | 62 | 11 | 18 | 54 | 47 | 26 | 3 | 27 | 2 | 55 | 46 |
| 21 | 48 | 57 | 4 | 60 | 1 | 24 | 45 | 41 | 20 | 5 | 64 | 8 | 61 | 44 | 17 |
| 40 | 29 | 12 | 49 | 9 | 52 | 37 | 32 | 28 | 33 | 56 | 13 | 53 | 16 | 25 | 36 |
| 59 | 34 | 23 | 14 | 22 | 15 | 58 | 35 | 7 | 30 | 43 | 50 | 42 | 51 | 6 | 31 |

132. Francis Gaspalou 2009 — D1= 98  D2= 162  D3= 162  D4= 98

| | | | | | | | | | | | | | | | |
|---|---|---|---|---|---|---|---|---|---|---|---|---|---|---|---|
| 10 | 19 | 38 | 63 | 39 | 62 | 11 | 18 | 22 | 15 | 58 | 35 | 59 | 34 | 23 | 14 |
| 37 | 64 | 9 | 20 | 12 | 17 | 40 | 61 | 57 | 36 | 21 | 16 | 24 | 13 | 60 | 33 |
| 56 | 45 | 28 | 1 | 25 | 4 | 53 | 48 | 44 | 49 | 8 | 29 | 5 | 32 | 41 | 52 |
| 27 | 2 | 55 | 46 | 54 | 47 | 26 | 3 | 7 | 30 | 43 | 50 | 42 | 51 | 6 | 31 |

133. Awani Kumar 2007 — D1= 66  D2= 194  D3= 66  D4= 194

| | | | | | | | | | | | | | | | |
|---|---|---|---|---|---|---|---|---|---|---|---|---|---|---|---|
| 10 | 19 | 38 | 63 | 57 | 48 | 21 | 4 | 40 | 29 | 12 | 49 | 23 | 34 | 59 | 14 |
| 39 | 62 | 11 | 18 | 24 | 1 | 60 | 45 | 9 | 52 | 37 | 32 | 58 | 15 | 22 | 35 |
| 26 | 47 | 54 | 3 | 41 | 20 | 5 | 64 | 56 | 33 | 28 | 13 | 7 | 30 | 43 | 50 |
| 55 | 2 | 27 | 46 | 8 | 61 | 44 | 17 | 25 | 16 | 53 | 36 | 42 | 51 | 6 | 31 |

134. Francis Gaspalou 2009 — D1= 70  D2= 198  D3= 126  D4= 126

| 10 | 19 | 38 | 63 |
|----|----|----|----|
| 39 | 62 | 11 | 18 |
| 54 | 47 | 26 | 3  |
| 27 | 2  | 55 | 46 |

| 53 | 16 | 25 | 36 |
|----|----|----|----|
| 28 | 33 | 56 | 13 |
| 9  | 52 | 37 | 32 |
| 40 | 29 | 12 | 49 |

| 8  | 61 | 44 | 17 |
|----|----|----|----|
| 41 | 20 | 5  | 64 |
| 60 | 1  | 24 | 45 |
| 21 | 48 | 57 | 4  |

| 59 | 34 | 23 | 14 |
|----|----|----|----|
| 22 | 15 | 58 | 35 |
| 7  | 30 | 43 | 50 |
| 42 | 51 | 6  | 31 |

135. Francis Gaspalou 2009    D1= 98    D2= 162    D3= 98    D4= 162

| 10 | 19 | 40 | 61 |
|----|----|----|----|
| 37 | 64 | 27 | 2  |
| 52 | 9  | 30 | 39 |
| 31 | 38 | 33 | 28 |

| 53 | 48 | 11 | 18 |
|----|----|----|----|
| 26 | 3  | 56 | 45 |
| 15 | 54 | 49 | 12 |
| 36 | 25 | 14 | 55 |

| 20 | 41 | 62 | 7  |
|----|----|----|----|
| 63 | 6  | 1  | 60 |
| 42 | 51 | 8  | 29 |
| 5  | 32 | 59 | 34 |

| 47 | 22 | 17 | 44 |
|----|----|----|----|
| 4  | 57 | 46 | 23 |
| 21 | 16 | 43 | 50 |
| 58 | 35 | 24 | 13 |

136. Awani Kumar 2017    D1= 34    D2= 226    D3= 130    D4= 130

| 10 | 19 | 54 | 47 |
|----|----|----|----|
| 37 | 64 | 9  | 20 |
| 56 | 45 | 28 | 1  |
| 27 | 2  | 39 | 62 |

| 55 | 46 | 11 | 18 |
|----|----|----|----|
| 12 | 17 | 40 | 61 |
| 25 | 4  | 53 | 48 |
| 38 | 63 | 26 | 3  |

| 6  | 15 | 58 | 51 |
|----|----|----|----|
| 57 | 36 | 21 | 16 |
| 44 | 49 | 8  | 29 |
| 23 | 30 | 43 | 34 |

| 59 | 50 | 7  | 14 |
|----|----|----|----|
| 24 | 13 | 60 | 33 |
| 5  | 32 | 41 | 52 |
| 42 | 35 | 22 | 31 |

137. Francis Gaspalou 2009    D1= 66    D2= 178    D3= 66    D4= 210

| 10 | 19 | 54 | 47 |
|----|----|----|----|
| 37 | 64 | 9  | 20 |
| 56 | 45 | 28 | 1  |
| 27 | 2  | 39 | 62 |

| 55 | 46 | 11 | 18 |
|----|----|----|----|
| 12 | 17 | 40 | 61 |
| 25 | 4  | 53 | 48 |
| 38 | 63 | 26 | 3  |

| 6  | 31 | 58 | 35 |
|----|----|----|----|
| 57 | 36 | 21 | 16 |
| 44 | 49 | 8  | 29 |
| 23 | 14 | 43 | 50 |

| 59 | 34 | 7  | 30 |
|----|----|----|----|
| 24 | 13 | 60 | 33 |
| 5  | 32 | 41 | 52 |
| 42 | 51 | 22 | 15 |

138. Francis Gaspalou 2009    D1= 50    D2= 178    D3= 82    D4= 210

| 10 | 19 | 54 | 47 |
|----|----|----|----|
| 37 | 64 | 9  | 20 |
| 56 | 45 | 28 | 1  |
| 27 | 2  | 39 | 62 |

| 55 | 46 | 11 | 18 |
|----|----|----|----|
| 12 | 17 | 40 | 61 |
| 25 | 4  | 53 | 48 |
| 38 | 63 | 26 | 3  |

| 22 | 15 | 42 | 51 |
|----|----|----|----|
| 57 | 36 | 21 | 16 |
| 44 | 49 | 8  | 29 |
| 7  | 30 | 59 | 34 |

| 43 | 50 | 23 | 14 |
|----|----|----|----|
| 24 | 13 | 60 | 33 |
| 5  | 32 | 41 | 52 |
| 58 | 35 | 6  | 31 |

139. Francis Gaspalou 2009    D1= 66    D2= 194    D3= 66    D4= 194

| 10 | 19 | 54 | 47 |
|----|----|----|----|
| 37 | 64 | 9  | 20 |
| 56 | 45 | 28 | 1  |
| 27 | 2  | 39 | 62 |

| 55 | 46 | 11 | 18 |
|----|----|----|----|
| 12 | 17 | 40 | 61 |
| 25 | 4  | 53 | 48 |
| 38 | 63 | 26 | 3  |

| 22 | 31 | 42 | 35 |
|----|----|----|----|
| 57 | 36 | 21 | 16 |
| 44 | 49 | 8  | 29 |
| 7  | 14 | 59 | 50 |

| 43 | 34 | 23 | 30 |
|----|----|----|----|
| 24 | 13 | 60 | 33 |
| 5  | 32 | 41 | 52 |
| 58 | 51 | 6  | 15 |

140. Francis Gaspalou 2009    D1= 50    D2= 194    D3= 82    D4= 194

| 10 | 19 | 56 | 45 |
|----|----|----|----|
| 53 | 48 | 11 | 18 |
| 4  | 57 | 62 | 7  |
| 63 | 6  | 1  | 60 |

| 55 | 46 | 9  | 20 |
|----|----|----|----|
| 12 | 17 | 54 | 47 |
| 29 | 40 | 35 | 26 |
| 34 | 27 | 32 | 37 |

| 52 | 41 | 14 | 23 |
|----|----|----|----|
| 15 | 22 | 49 | 44 |
| 58 | 3  | 8  | 61 |
| 5  | 64 | 59 | 2  |

| 13 | 24 | 51 | 42 |
|----|----|----|----|
| 50 | 43 | 16 | 21 |
| 39 | 30 | 25 | 36 |
| 28 | 33 | 38 | 31 |

141. Awani Kumar 2017    D1= 66    D2= 130    D3= 194    D4= 130

| 10 | 21 | 40 | 59 |
|----|----|----|----|
| 39 | 60 | 9  | 22 |
| 54 | 41 | 28 | 7  |
| 27 | 8  | 53 | 42 |

| 51 | 48 | 29 | 2  |
|----|----|----|----|
| 30 | 1  | 52 | 47 |
| 15 | 20 | 33 | 62 |
| 34 | 61 | 16 | 19 |

| 38 | 57 | 12 | 23 |
|----|----|----|----|
| 11 | 24 | 37 | 58 |
| 26 | 5  | 56 | 43 |
| 55 | 44 | 25 | 6  |

| 31 | 4  | 49 | 46 |
|----|----|----|----|
| 50 | 45 | 32 | 3  |
| 35 | 64 | 13 | 18 |
| 14 | 17 | 36 | 63 |

142. Awani Kumar 2007    D1= 130    D2= 130    D3= 130    D4= 130

| 10 | 21 | 40 | 59 |
|----|----|----|----|
| 23 | 60 | 9  | 38 |
| 54 | 41 | 28 | 7  |
| 43 | 8  | 53 | 26 |

| 51 | 48 | 29 | 2  |
|----|----|----|----|
| 46 | 1  | 52 | 31 |
| 15 | 20 | 33 | 62 |
| 18 | 61 | 16 | 35 |

| 22 | 57 | 12 | 39 |
|----|----|----|----|
| 11 | 24 | 37 | 58 |
| 42 | 5  | 56 | 27 |
| 55 | 44 | 25 | 6  |

| 47 | 4  | 49 | 30 |
|----|----|----|----|
| 50 | 45 | 32 | 3  |
| 19 | 64 | 13 | 34 |
| 14 | 17 | 36 | 63 |

143. Francis Gaspalou 2009    D1= 130    D2= 130    D3= 130    D4= 130

| 10 | 21 | 40 | 59 | 63 | 36 | 17 | 14 | 4 | 57 | 46 | 23 | 53 | 16 | 27 | 34 |
|---|---|---|---|---|---|---|---|---|---|---|---|---|---|---|---|
| 39 | 60 | 9 | 22 | 18 | 13 | 64 | 35 | 45 | 24 | 3 | 58 | 28 | 33 | 54 | 15 |
| 62 | 41 | 20 | 7 | 11 | 32 | 37 | 50 | 56 | 5 | 26 | 43 | 1 | 52 | 47 | 30 |
| 19 | 8 | 61 | 42 | 38 | 49 | 12 | 31 | 25 | 44 | 55 | 6 | 48 | 29 | 2 | 51 |

144. Awani Kumar 2011    D1= 100    D2= 176    D3= 88    D4= 156

| 10 | 29 | 52 | 39 | 53 | 34 | 15 | 28 | 4 | 23 | 58 | 45 | 63 | 44 | 5 | 18 |
|---|---|---|---|---|---|---|---|---|---|---|---|---|---|---|---|
| 51 | 40 | 9 | 30 | 16 | 27 | 54 | 33 | 57 | 46 | 3 | 24 | 6 | 17 | 64 | 43 |
| 56 | 35 | 14 | 25 | 11 | 32 | 49 | 38 | 62 | 41 | 8 | 19 | 1 | 22 | 59 | 48 |
| 13 | 26 | 55 | 36 | 50 | 37 | 12 | 31 | 7 | 20 | 61 | 42 | 60 | 47 | 2 | 21 |

145. Francis Gaspalou 2009    D1= 66    D2= 194    D3= 66    D4= 194

| 10 | 31 | 40 | 49 | 37 | 52 | 11 | 30 | 32 | 41 | 50 | 7 | 51 | 6 | 29 | 44 |
|---|---|---|---|---|---|---|---|---|---|---|---|---|---|---|---|
| 33 | 56 | 15 | 26 | 14 | 27 | 36 | 53 | 55 | 2 | 25 | 48 | 28 | 45 | 54 | 3 |
| 64 | 9 | 18 | 39 | 19 | 38 | 61 | 12 | 42 | 63 | 8 | 17 | 5 | 20 | 43 | 62 |
| 23 | 34 | 57 | 16 | 60 | 13 | 22 | 35 | 1 | 24 | 47 | 58 | 46 | 59 | 4 | 21 |

146. Guenter Stertenbrink 2003    D1= 66    D2= 194    D3= 130    D4= 130

| 10 | 33 | 64 | 23 | 63 | 24 | 41 | 2 | 40 | 15 | 18 | 57 | 17 | 58 | 7 | 48 |
|---|---|---|---|---|---|---|---|---|---|---|---|---|---|---|---|
| 37 | 14 | 19 | 60 | 20 | 59 | 6 | 45 | 11 | 36 | 61 | 22 | 62 | 21 | 44 | 3 |
| 32 | 55 | 42 | 1 | 9 | 34 | 31 | 56 | 50 | 25 | 8 | 47 | 39 | 16 | 49 | 26 |
| 51 | 28 | 5 | 46 | 38 | 13 | 52 | 27 | 29 | 54 | 43 | 4 | 12 | 35 | 30 | 53 |

147. Guenter Stertenbrink 2003    D1= 130    D2= 66    D3= 194    D4= 130

| 10 | 35 | 54 | 31 | 55 | 30 | 11 | 34 | 38 | 15 | 26 | 51 | 27 | 50 | 39 | 14 |
|---|---|---|---|---|---|---|---|---|---|---|---|---|---|---|---|
| 53 | 16 | 25 | 36 | 28 | 33 | 56 | 13 | 9 | 52 | 37 | 32 | 40 | 29 | 12 | 49 |
| 8 | 61 | 44 | 17 | 41 | 20 | 5 | 64 | 60 | 1 | 24 | 45 | 21 | 48 | 57 | 4 |
| 59 | 18 | 7 | 46 | 6 | 47 | 58 | 19 | 23 | 62 | 43 | 2 | 42 | 3 | 22 | 63 |

148. Awani Kumar 2009    D1= 130    D2= 130    D3= 130    D4= 130

| 10 | 35 | 54 | 31 | 55 | 30 | 11 | 34 | 38 | 15 | 26 | 51 | 27 | 50 | 39 | 14 |
|---|---|---|---|---|---|---|---|---|---|---|---|---|---|---|---|
| 53 | 16 | 25 | 36 | 28 | 33 | 56 | 13 | 9 | 52 | 37 | 32 | 40 | 29 | 12 | 49 |
| 24 | 61 | 44 | 1 | 41 | 4 | 21 | 64 | 60 | 17 | 8 | 45 | 5 | 48 | 57 | 20 |
| 43 | 18 | 7 | 62 | 6 | 63 | 42 | 19 | 23 | 46 | 59 | 2 | 58 | 3 | 22 | 47 |

149. Awani Kumar 2007    D1= 98    D2= 162    D3= 98    D4= 162

| 10 | 37 | 52 | 31 | 53 | 26 | 15 | 36 | 20 | 63 | 42 | 5 | 47 | 4 | 21 | 58 |
|---|---|---|---|---|---|---|---|---|---|---|---|---|---|---|---|
| 51 | 32 | 41 | 6 | 16 | 35 | 22 | 57 | 9 | 38 | 19 | 64 | 54 | 25 | 48 | 3 |
| 40 | 27 | 30 | 33 | 11 | 56 | 49 | 14 | 62 | 1 | 8 | 59 | 17 | 46 | 43 | 24 |
| 29 | 34 | 7 | 60 | 50 | 13 | 44 | 23 | 39 | 28 | 61 | 2 | 12 | 55 | 18 | 45 |

150. Awani Kumar 2017    D1= 98    D2= 66    D3= 162    D4= 194

| 10 | 39 | 26 | 55 | 57 | 24 | 41 | 8 | 40 | 9 | 56 | 25 | 23 | 58 | 7 | 42 |
|---|---|---|---|---|---|---|---|---|---|---|---|---|---|---|---|
| 51 | 30 | 15 | 34 | 48 | 1 | 20 | 61 | 29 | 52 | 33 | 16 | 2 | 47 | 62 | 19 |
| 38 | 11 | 54 | 27 | 21 | 60 | 5 | 44 | 12 | 37 | 28 | 53 | 59 | 22 | 43 | 6 |
| 31 | 50 | 35 | 14 | 4 | 45 | 64 | 17 | 49 | 32 | 13 | 36 | 46 | 3 | 18 | 63 |

151. Awani Kumar 2007    D1= 102    D2= 158    D3= 166    D4= 94

| 10 | 39 | 52 | 29 | 53 | 28 | 15 | 34 | 20 | 45 | 42 | 23 | 47 | 18 | 21 | 44 |
|---|---|---|---|---|---|---|---|---|---|---|---|---|---|---|---|
| 51 | 30 | 9 | 40 | 16 | 33 | 54 | 27 | 57 | 8 | 3 | 62 | 6 | 59 | 64 | 1 |
| 56 | 25 | 14 | 35 | 11 | 38 | 49 | 32 | 46 | 19 | 24 | 41 | 17 | 48 | 43 | 22 |
| 13 | 36 | 55 | 26 | 50 | 31 | 12 | 37 | 7 | 58 | 61 | 4 | 60 | 5 | 2 | 63 |

152. Awani Kumar 2012    D1= 130    D2= 162    D3= 98    D4= 130

| | | | | | | | | | | | | | | | |
|---|---|---|---|---|---|---|---|---|---|---|---|---|---|---|---|
| 10 | 39 | 54 | 27 | 53 | 28 | 9 | 40 | 8 | 41 | 60 | 21 | 59 | 22 | 7 | 42 |
| 51 | 30 | 15 | 34 | 16 | 33 | 52 | 29 | 61 | 20 | 1 | 48 | 2 | 47 | 62 | 19 |
| 38 | 11 | 26 | 55 | 25 | 56 | 37 | 12 | 44 | 5 | 24 | 57 | 23 | 58 | 43 | 6 |
| 31 | 50 | 35 | 14 | 36 | 13 | 32 | 49 | 17 | 64 | 45 | 4 | 46 | 3 | 18 | 63 |

153. Francis Gaspalou 2009    D1= 130    D2= 130    D3= 130    D4= 130

| | | | | | | | | | | | | | | | |
|---|---|---|---|---|---|---|---|---|---|---|---|---|---|---|---|
| 10 | 39 | 54 | 27 | 53 | 28 | 9 | 40 | 8 | 57 | 44 | 21 | 59 | 6 | 23 | 42 |
| 51 | 30 | 15 | 34 | 16 | 33 | 52 | 29 | 45 | 20 | 1 | 64 | 18 | 47 | 62 | 3 |
| 38 | 11 | 26 | 55 | 25 | 56 | 37 | 12 | 60 | 5 | 24 | 41 | 7 | 58 | 43 | 22 |
| 31 | 50 | 35 | 14 | 36 | 13 | 32 | 49 | 17 | 48 | 61 | 4 | 46 | 19 | 2 | 63 |

154. Awani Kumar 2009    D1= 130    D2= 130    D3= 130    D4= 130

| | | | | | | | | | | | | | | | |
|---|---|---|---|---|---|---|---|---|---|---|---|---|---|---|---|
| 10 | 39 | 54 | 27 | 53 | 28 | 9 | 40 | 24 | 41 | 60 | 5 | 43 | 22 | 7 | 58 |
| 51 | 30 | 15 | 34 | 16 | 33 | 52 | 29 | 61 | 4 | 17 | 48 | 2 | 63 | 46 | 19 |
| 38 | 11 | 26 | 55 | 25 | 56 | 37 | 12 | 44 | 21 | 8 | 57 | 23 | 42 | 59 | 6 |
| 31 | 50 | 35 | 14 | 36 | 13 | 32 | 49 | 1 | 64 | 45 | 20 | 62 | 3 | 18 | 47 |

155. Awani Kumar 2009    D1= 98    D2= 162    D3= 162    D4= 98

| | | | | | | | | | | | | | | | |
|---|---|---|---|---|---|---|---|---|---|---|---|---|---|---|---|
| 10 | 39 | 54 | 27 | 53 | 28 | 9 | 40 | 24 | 57 | 44 | 5 | 43 | 6 | 23 | 58 |
| 51 | 30 | 15 | 34 | 16 | 33 | 52 | 29 | 45 | 4 | 17 | 64 | 18 | 63 | 46 | 3 |
| 38 | 11 | 26 | 55 | 25 | 56 | 37 | 12 | 60 | 21 | 8 | 41 | 7 | 42 | 59 | 22 |
| 31 | 50 | 35 | 14 | 36 | 13 | 32 | 49 | 1 | 48 | 61 | 20 | 62 | 19 | 2 | 47 |

156. Awani Kumar 2007    D1= 98    D2= 162    D3= 162    D4= 98

| | | | | | | | | | | | | | | | |
|---|---|---|---|---|---|---|---|---|---|---|---|---|---|---|---|
| 10 | 51 | 54 | 15 | 55 | 14 | 11 | 50 | 38 | 31 | 26 | 35 | 27 | 34 | 39 | 30 |
| 53 | 16 | 25 | 36 | 28 | 33 | 56 | 13 | 9 | 52 | 37 | 32 | 40 | 29 | 12 | 49 |
| 8 | 61 | 44 | 17 | 41 | 20 | 5 | 64 | 60 | 1 | 24 | 45 | 21 | 48 | 57 | 4 |
| 59 | 2 | 7 | 62 | 6 | 63 | 58 | 3 | 23 | 46 | 43 | 18 | 42 | 19 | 22 | 47 |

157. Awani Kumar 2009    D1= 114    D2= 114    D3= 146    D4= 146

| | | | | | | | | | | | | | | | |
|---|---|---|---|---|---|---|---|---|---|---|---|---|---|---|---|
| 11 | 2 | 55 | 62 | 56 | 45 | 28 | 1 | 37 | 64 | 9 | 20 | 26 | 19 | 38 | 47 |
| 54 | 63 | 10 | 3 | 25 | 4 | 53 | 48 | 12 | 17 | 40 | 61 | 39 | 46 | 27 | 18 |
| 7 | 14 | 59 | 50 | 44 | 49 | 8 | 29 | 57 | 36 | 21 | 16 | 22 | 31 | 42 | 35 |
| 58 | 51 | 6 | 15 | 5 | 32 | 41 | 52 | 24 | 13 | 60 | 33 | 43 | 34 | 23 | 30 |

158. Awani Kumar 2009    D1= 66    D2= 194    D3= 194    D4= 66

| | | | | | | | | | | | | | | | |
|---|---|---|---|---|---|---|---|---|---|---|---|---|---|---|---|
| 11 | 2 | 55 | 62 | 56 | 45 | 28 | 1 | 37 | 64 | 9 | 20 | 26 | 19 | 38 | 47 |
| 54 | 63 | 10 | 3 | 25 | 4 | 53 | 48 | 12 | 17 | 40 | 61 | 39 | 46 | 27 | 18 |
| 7 | 30 | 59 | 34 | 44 | 49 | 8 | 29 | 57 | 36 | 21 | 16 | 22 | 15 | 42 | 51 |
| 58 | 35 | 6 | 31 | 5 | 32 | 41 | 52 | 24 | 13 | 60 | 33 | 43 | 50 | 23 | 14 |

159. Francis Gaspalou 2009    D1= 50    D2= 194    D3= 194    D4= 82

| | | | | | | | | | | | | | | | |
|---|---|---|---|---|---|---|---|---|---|---|---|---|---|---|---|
| 11 | 2 | 55 | 62 | 56 | 45 | 28 | 1 | 37 | 64 | 9 | 20 | 26 | 19 | 38 | 47 |
| 54 | 63 | 10 | 3 | 25 | 4 | 53 | 48 | 12 | 17 | 40 | 61 | 39 | 46 | 27 | 18 |
| 23 | 14 | 43 | 50 | 44 | 49 | 8 | 29 | 57 | 36 | 21 | 16 | 6 | 31 | 58 | 35 |
| 42 | 51 | 22 | 15 | 5 | 32 | 41 | 52 | 24 | 13 | 60 | 33 | 59 | 34 | 7 | 30 |

160. Francis Gaspalou 2009    D1= 66    D2= 210    D3= 178    D4= 66

| | | | | | | | | | | | | | | | |
|---|---|---|---|---|---|---|---|---|---|---|---|---|---|---|---|
| 11 | 2 | 55 | 62 | 56 | 45 | 28 | 1 | 37 | 64 | 9 | 20 | 26 | 19 | 38 | 47 |
| 54 | 63 | 10 | 3 | 25 | 4 | 53 | 48 | 12 | 17 | 40 | 61 | 39 | 46 | 27 | 18 |
| 23 | 30 | 43 | 34 | 44 | 49 | 8 | 29 | 57 | 36 | 21 | 16 | 6 | 15 | 58 | 51 |
| 42 | 35 | 22 | 31 | 5 | 32 | 41 | 52 | 24 | 13 | 60 | 33 | 59 | 50 | 7 | 14 |

161. Francis Gaspalou 2009    D1= 50    D2= 210    D3= 178    D4= 82

| | | | | | | | | | | | | | | | |
|---|---|---|---|---|---|---|---|---|---|---|---|---|---|---|---|
| 11 | 18 | 39 | 62 | 38 | 63 | 10 | 19 | 55 | 46 | 27 | 2 | 26 | 3 | 54 | 47 |
| 24 | 45 | 28 | 33 | 57 | 4 | 53 | 16 | 44 | 17 | 40 | 29 | 5 | 64 | 9 | 52 |
| 37 | 32 | 41 | 20 | 12 | 49 | 8 | 61 | 25 | 36 | 21 | 48 | 56 | 13 | 60 | 1 |
| 58 | 35 | 22 | 15 | 23 | 14 | 59 | 34 | 6 | 31 | 42 | 51 | 43 | 50 | 7 | 30 |
| 162. Awani Kumar 2017 | | | | D1= | 66 | D2= | 194 | | D3= | 194 | | D4= | 66 | | |
| 11 | 18 | 39 | 62 | 56 | 45 | 28 | 1 | 37 | 64 | 9 | 20 | 26 | 3 | 54 | 47 |
| 38 | 63 | 10 | 19 | 25 | 4 | 53 | 48 | 12 | 17 | 40 | 61 | 55 | 46 | 27 | 2 |
| 23 | 14 | 59 | 34 | 44 | 49 | 8 | 29 | 57 | 36 | 21 | 16 | 6 | 31 | 42 | 51 |
| 58 | 35 | 22 | 15 | 5 | 32 | 41 | 52 | 24 | 13 | 60 | 33 | 43 | 50 | 7 | 30 |
| 163. Awani Kumar 2007 | | | | D1= | 66 | D2= | 194 | | D3= | 194 | | D4= | 66 | | |
| 11 | 18 | 39 | 62 | 56 | 13 | 60 | 1 | 5 | 64 | 9 | 52 | 58 | 35 | 22 | 15 |
| 38 | 63 | 10 | 19 | 25 | 36 | 21 | 48 | 44 | 17 | 40 | 29 | 23 | 14 | 59 | 34 |
| 55 | 46 | 27 | 2 | 12 | 49 | 8 | 61 | 57 | 4 | 53 | 16 | 6 | 31 | 42 | 51 |
| 26 | 3 | 54 | 47 | 37 | 32 | 41 | 20 | 24 | 45 | 28 | 33 | 43 | 50 | 7 | 30 |
| 164. Awani Kumar 2017 | | | | D1= | 130 | D2= | 130 | | D3= | 130 | | D4= | 130 | | |
| 11 | 18 | 55 | 46 | 56 | 45 | 28 | 1 | 37 | 64 | 9 | 20 | 26 | 3 | 38 | 63 |
| 54 | 47 | 10 | 19 | 25 | 4 | 53 | 48 | 12 | 17 | 40 | 61 | 39 | 62 | 27 | 2 |
| 7 | 14 | 59 | 50 | 44 | 49 | 8 | 29 | 57 | 36 | 21 | 16 | 22 | 31 | 42 | 35 |
| 58 | 51 | 6 | 15 | 5 | 32 | 41 | 52 | 24 | 13 | 60 | 33 | 43 | 34 | 23 | 30 |
| 165. Francis Gaspalou 2009 | | | | D1= | 66 | D2= | 178 | | D3= | 210 | | D4= | 66 | | |
| 11 | 18 | 55 | 46 | 56 | 45 | 28 | 1 | 37 | 64 | 9 | 20 | 26 | 3 | 38 | 63 |
| 54 | 47 | 10 | 19 | 25 | 4 | 53 | 48 | 12 | 17 | 40 | 61 | 39 | 62 | 27 | 2 |
| 7 | 30 | 59 | 34 | 44 | 49 | 8 | 29 | 57 | 36 | 21 | 16 | 22 | 15 | 42 | 51 |
| 58 | 35 | 6 | 31 | 5 | 32 | 41 | 52 | 24 | 13 | 60 | 33 | 43 | 50 | 23 | 14 |
| 166. Francis Gaspalou 2009 | | | | D1= | 50 | D2= | 178 | | D3= | 210 | | D4= | 82 | | |
| 11 | 18 | 55 | 46 | 56 | 45 | 28 | 1 | 37 | 64 | 9 | 20 | 26 | 3 | 38 | 63 |
| 54 | 47 | 10 | 19 | 25 | 4 | 53 | 48 | 12 | 17 | 40 | 61 | 39 | 62 | 27 | 2 |
| 23 | 14 | 43 | 50 | 44 | 49 | 8 | 29 | 57 | 36 | 21 | 16 | 6 | 31 | 58 | 35 |
| 42 | 51 | 22 | 15 | 5 | 32 | 41 | 52 | 24 | 13 | 60 | 33 | 59 | 34 | 7 | 30 |
| 167. Francis Gaspalou 2009 | | | | D1= | 66 | D2= | 194 | | D3= | 194 | | D4= | 66 | | |
| 11 | 18 | 55 | 46 | 56 | 45 | 28 | 1 | 37 | 64 | 9 | 20 | 26 | 3 | 38 | 63 |
| 54 | 47 | 10 | 19 | 25 | 4 | 53 | 48 | 12 | 17 | 40 | 61 | 39 | 62 | 27 | 2 |
| 23 | 30 | 43 | 34 | 44 | 49 | 8 | 29 | 57 | 36 | 21 | 16 | 6 | 15 | 58 | 51 |
| 42 | 35 | 22 | 31 | 5 | 32 | 41 | 52 | 24 | 13 | 60 | 33 | 59 | 50 | 7 | 14 |
| 168. Francis Gaspalou 2009 | | | | D1= | 50 | D2= | 194 | | D3= | 194 | | D4= | 82 | | |
| 11 | 22 | 55 | 42 | 50 | 47 | 14 | 19 | 23 | 10 | 43 | 54 | 46 | 51 | 18 | 15 |
| 40 | 9 | 28 | 53 | 29 | 52 | 33 | 16 | 12 | 37 | 56 | 25 | 49 | 32 | 13 | 36 |
| 21 | 60 | 41 | 8 | 48 | 1 | 20 | 61 | 57 | 24 | 5 | 44 | 4 | 45 | 64 | 17 |
| 58 | 39 | 6 | 27 | 3 | 30 | 63 | 34 | 38 | 59 | 26 | 7 | 31 | 2 | 35 | 62 |
| 169. Francis Gaspalou 2009 | | | | D1= | 130 | D2= | 130 | | D3= | 130 | | D4= | 130 | | |
| 11 | 24 | 37 | 58 | 50 | 45 | 32 | 3 | 39 | 28 | 41 | 22 | 30 | 33 | 20 | 47 |
| 38 | 57 | 12 | 23 | 31 | 4 | 49 | 46 | 10 | 53 | 8 | 59 | 51 | 16 | 61 | 2 |
| 55 | 44 | 25 | 6 | 14 | 17 | 36 | 63 | 27 | 40 | 21 | 42 | 34 | 29 | 48 | 19 |
| 26 | 5 | 56 | 43 | 35 | 64 | 13 | 18 | 54 | 9 | 60 | 7 | 15 | 52 | 1 | 62 |
| 170. Awani Kumar 2012 | | | | D1= | 98 | D2= | 162 | | D3= | 98 | | D4= | 162 | | |

| 11 | 34 | 55 | 30 | | 54 | 31 | 10 | 35 | | 39 | 14 | 27 | 50 | | 26 | 51 | 38 | 15 |
|---|---|---|---|---|---|---|---|---|---|---|---|---|---|---|---|---|---|---|
| 40 | 29 | 12 | 49 | | 9 | 52 | 37 | 32 | | 28 | 33 | 56 | 13 | | 53 | 16 | 25 | 36 |
| 21 | 48 | 57 | 4 | | 60 | 1 | 24 | 45 | | 41 | 20 | 5 | 64 | | 8 | 61 | 44 | 17 |
| 58 | 19 | 6 | 47 | | 7 | 46 | 59 | 18 | | 22 | 63 | 42 | 3 | | 43 | 2 | 23 | 62 |

171. Francis Gaspalou 2009    D1= 130    D2= 130    D3= 130    D4= 130

| 11 | 34 | 55 | 30 | | 54 | 31 | 10 | 35 | | 39 | 14 | 27 | 50 | | 26 | 51 | 38 | 15 |
|---|---|---|---|---|---|---|---|---|---|---|---|---|---|---|---|---|---|---|
| 40 | 29 | 12 | 49 | | 9 | 52 | 37 | 32 | | 28 | 33 | 56 | 13 | | 53 | 16 | 25 | 36 |
| 21 | 64 | 41 | 4 | | 44 | 1 | 24 | 61 | | 57 | 20 | 5 | 48 | | 8 | 45 | 60 | 17 |
| 58 | 3 | 22 | 47 | | 23 | 46 | 59 | 2 | | 6 | 63 | 42 | 19 | | 43 | 18 | 7 | 62 |

172. Awani Kumar 2007    D1= 130    D2= 130    D3= 130    D4= 130

| 11 | 38 | 27 | 54 | | 50 | 31 | 14 | 35 | | 39 | 10 | 55 | 26 | | 30 | 51 | 34 | 15 |
|---|---|---|---|---|---|---|---|---|---|---|---|---|---|---|---|---|---|---|
| 40 | 9 | 56 | 25 | | 29 | 52 | 33 | 16 | | 12 | 37 | 28 | 53 | | 49 | 32 | 13 | 36 |
| 57 | 24 | 41 | 8 | | 48 | 1 | 20 | 61 | | 21 | 60 | 5 | 44 | | 4 | 45 | 64 | 17 |
| 22 | 59 | 6 | 43 | | 3 | 46 | 63 | 18 | | 58 | 23 | 42 | 7 | | 47 | 2 | 19 | 62 |

173. Awani Kumar 2007    D1= 130    D2= 194    D3= 66    D4= 130

| 11 | 38 | 55 | 26 | | 50 | 31 | 14 | 35 | | 39 | 10 | 27 | 54 | | 30 | 51 | 34 | 15 |
|---|---|---|---|---|---|---|---|---|---|---|---|---|---|---|---|---|---|---|
| 40 | 9 | 28 | 53 | | 29 | 52 | 33 | 16 | | 12 | 37 | 56 | 25 | | 49 | 32 | 13 | 36 |
| 21 | 60 | 41 | 8 | | 48 | 1 | 20 | 61 | | 57 | 24 | 5 | 44 | | 4 | 45 | 64 | 17 |
| 58 | 23 | 6 | 43 | | 3 | 46 | 63 | 18 | | 22 | 59 | 42 | 7 | | 47 | 2 | 19 | 62 |

174. Awani Kumar 2007    D1= 130    D2= 130    D3= 130    D4= 130

| 11 | 38 | 55 | 26 | | 56 | 25 | 12 | 37 | | 5 | 44 | 57 | 24 | | 58 | 23 | 6 | 43 |
|---|---|---|---|---|---|---|---|---|---|---|---|---|---|---|---|---|---|---|
| 50 | 31 | 14 | 35 | | 13 | 36 | 49 | 32 | | 64 | 17 | 4 | 45 | | 3 | 46 | 63 | 18 |
| 39 | 10 | 27 | 54 | | 60 | 21 | 8 | 41 | | 9 | 40 | 53 | 28 | | 22 | 59 | 42 | 7 |
| 30 | 51 | 34 | 15 | | 1 | 48 | 61 | 20 | | 52 | 29 | 16 | 33 | | 47 | 2 | 19 | 62 |

175. Awani Kumar 2012    D1= 162    D2= 162    D3= 98    D4= 98

| 11 | 40 | 21 | 58 | | 54 | 9 | 60 | 7 | | 39 | 28 | 41 | 22 | | 26 | 53 | 8 | 43 |
|---|---|---|---|---|---|---|---|---|---|---|---|---|---|---|---|---|---|---|
| 50 | 29 | 48 | 3 | | 15 | 52 | 1 | 62 | | 30 | 33 | 20 | 47 | | 35 | 16 | 61 | 18 |
| 55 | 12 | 57 | 6 | | 10 | 37 | 24 | 59 | | 27 | 56 | 5 | 42 | | 38 | 25 | 44 | 23 |
| 14 | 49 | 4 | 63 | | 51 | 32 | 45 | 2 | | 34 | 13 | 64 | 19 | | 31 | 36 | 17 | 46 |

176. Awani Kumar 2009    D1= 114    D2= 146    D3= 114    D4= 146

| 12 | 17 | 54 | 47 | | 55 | 46 | 9 | 20 | | 2 | 59 | 64 | 5 | | 61 | 8 | 3 | 58 |
|---|---|---|---|---|---|---|---|---|---|---|---|---|---|---|---|---|---|---|
| 53 | 48 | 11 | 18 | | 10 | 19 | 56 | 45 | | 31 | 38 | 33 | 28 | | 36 | 25 | 30 | 39 |
| 50 | 43 | 16 | 21 | | 13 | 24 | 51 | 42 | | 60 | 1 | 6 | 63 | | 7 | 62 | 57 | 4 |
| 15 | 22 | 49 | 44 | | 52 | 41 | 14 | 23 | | 37 | 32 | 27 | 34 | | 26 | 35 | 40 | 29 |

177. Awani Kumar 2017    D1= 66    D2= 130    D3= 130    D4= 194

| 12 | 31 | 50 | 37 | | 55 | 36 | 13 | 26 | | 2 | 21 | 60 | 47 | | 61 | 42 | 7 | 20 |
|---|---|---|---|---|---|---|---|---|---|---|---|---|---|---|---|---|---|---|
| 49 | 38 | 11 | 32 | | 14 | 25 | 56 | 35 | | 59 | 48 | 1 | 22 | | 8 | 19 | 62 | 41 |
| 54 | 33 | 16 | 27 | | 9 | 30 | 51 | 40 | | 64 | 43 | 6 | 17 | | 3 | 24 | 57 | 46 |
| 15 | 28 | 53 | 34 | | 52 | 39 | 10 | 29 | | 5 | 18 | 63 | 44 | | 58 | 45 | 4 | 23 |

178. Francis Gaspalou 2009    D1= 66    D2= 194    D3= 66    D4= 194

| 12 | 35 | 30 | 53 | | 61 | 22 | 11 | 36 | | 38 | 13 | 52 | 27 | | 19 | 60 | 37 | 14 |
|---|---|---|---|---|---|---|---|---|---|---|---|---|---|---|---|---|---|---|
| 39 | 16 | 49 | 26 | | 18 | 57 | 40 | 15 | | 9 | 34 | 31 | 56 | | 64 | 23 | 10 | 33 |
| 62 | 21 | 44 | 3 | | 43 | 4 | 29 | 54 | | 20 | 59 | 6 | 45 | | 5 | 46 | 51 | 28 |
| 17 | 58 | 7 | 48 | | 8 | 47 | 50 | 25 | | 63 | 24 | 41 | 2 | | 42 | 1 | 32 | 55 |

179. Guenter Stertenbrink 2003    D1= 130    D2= 194    D3= 66    D4= 130

| 12 | 37 | 50 | 31 | | 55 | 26 | 13 | 36 | | 18 | 47 | 44 | 21 | | 45 | 20 | 23 | 42 |
|---|---|---|---|---|---|---|---|---|---|---|---|---|---|---|---|---|---|---|
| 49 | 32 | 11 | 38 | | 14 | 35 | 56 | 25 | | 59 | 6 | 1 | 64 | | 8 | 57 | 62 | 3 |
| 54 | 27 | 16 | 33 | | 9 | 40 | 51 | 30 | | 48 | 17 | 22 | 43 | | 19 | 46 | 41 | 24 |
| 15 | 34 | 53 | 28 | | 52 | 29 | 10 | 39 | | 5 | 60 | 63 | 2 | | 58 | 7 | 4 | 61 |

180. Awani Kumar 2012    D1= 130    D2= 162    D3= 98    D4= 130

| 12 | 39 | 50 | 29 | | 55 | 28 | 13 | 34 | | 18 | 61 | 44 | 7 | | 45 | 2 | 23 | 60 |
|---|---|---|---|---|---|---|---|---|---|---|---|---|---|---|---|---|---|---|
| 49 | 30 | 43 | 8 | | 14 | 33 | 24 | 59 | | 11 | 40 | 17 | 62 | | 56 | 27 | 46 | 1 |
| 54 | 9 | 16 | 51 | | 25 | 38 | 35 | 32 | | 48 | 19 | 22 | 41 | | 3 | 64 | 57 | 6 |
| 15 | 52 | 21 | 42 | | 36 | 31 | 58 | 5 | | 53 | 10 | 47 | 20 | | 26 | 37 | 4 | 63 |

181. Awani Kumar 2017    D1= 130    D2= 98    D3= 130    D4= 162

| 13 | 8 | 51 | 58 | | 50 | 59 | 16 | 5 | | 7 | 14 | 57 | 52 | | 60 | 49 | 6 | 15 |
|---|---|---|---|---|---|---|---|---|---|---|---|---|---|---|---|---|---|---|
| 26 | 19 | 40 | 45 | | 37 | 48 | 27 | 18 | | 20 | 25 | 46 | 39 | | 47 | 38 | 17 | 28 |
| 55 | 62 | 9 | 4 | | 12 | 1 | 54 | 63 | | 61 | 56 | 3 | 10 | | 2 | 11 | 64 | 53 |
| 36 | 41 | 30 | 23 | | 31 | 22 | 33 | 44 | | 42 | 35 | 24 | 29 | | 21 | 32 | 43 | 34 |

182. Awani Kumar 2009    D1= 98    D2= 162    D3= 98    D4= 162

| 13 | 34 | 55 | 28 | | 56 | 27 | 46 | 1 | | 35 | 32 | 25 | 38 | | 26 | 37 | 4 | 63 |
|---|---|---|---|---|---|---|---|---|---|---|---|---|---|---|---|---|---|---|
| 50 | 29 | 12 | 39 | | 11 | 40 | 17 | 62 | | 16 | 51 | 54 | 9 | | 53 | 10 | 47 | 20 |
| 23 | 60 | 45 | 2 | | 14 | 33 | 24 | 59 | | 57 | 6 | 3 | 64 | | 36 | 31 | 58 | 5 |
| 44 | 7 | 18 | 61 | | 49 | 30 | 43 | 8 | | 22 | 41 | 48 | 19 | | 15 | 52 | 21 | 42 |

183. Awani Kumar 2017    D1= 98    D2= 66    D3= 194    D4= 162

| 14 | 33 | 28 | 55 | | 59 | 24 | 13 | 34 | | 36 | 15 | 54 | 25 | | 21 | 58 | 35 | 16 |
|---|---|---|---|---|---|---|---|---|---|---|---|---|---|---|---|---|---|---|
| 37 | 10 | 51 | 32 | | 20 | 63 | 38 | 9 | | 11 | 40 | 29 | 50 | | 62 | 17 | 12 | 39 |
| 60 | 23 | 46 | 1 | | 45 | 2 | 27 | 56 | | 22 | 57 | 4 | 47 | | 3 | 48 | 53 | 26 |
| 19 | 64 | 5 | 42 | | 6 | 41 | 52 | 31 | | 61 | 18 | 43 | 8 | | 44 | 7 | 30 | 49 |

184. Guenter Stertenbrink 2003    D1= 130    D2= 194    D3= 66    D4= 130

| 18 | 3 | 64 | 45 | | 63 | 46 | 17 | 4 | | 14 | 31 | 36 | 49 | | 35 | 50 | 13 | 32 |
|---|---|---|---|---|---|---|---|---|---|---|---|---|---|---|---|---|---|---|
| 11 | 26 | 37 | 56 | | 38 | 55 | 12 | 25 | | 23 | 6 | 57 | 44 | | 58 | 43 | 24 | 5 |
| 62 | 47 | 20 | 1 | | 19 | 2 | 61 | 48 | | 34 | 51 | 16 | 29 | | 15 | 30 | 33 | 52 |
| 39 | 54 | 9 | 28 | | 10 | 27 | 40 | 53 | | 59 | 42 | 21 | 8 | | 22 | 7 | 60 | 41 |

185. Awani Kumar 2006    D1= 130    D2= 130    D3= 130    D4= 130

| 18 | 3 | 64 | 45 | | 63 | 46 | 17 | 4 | | 14 | 31 | 36 | 49 | | 35 | 50 | 13 | 32 |
|---|---|---|---|---|---|---|---|---|---|---|---|---|---|---|---|---|---|---|
| 11 | 54 | 37 | 28 | | 38 | 27 | 12 | 53 | | 59 | 6 | 21 | 44 | | 22 | 43 | 60 | 5 |
| 62 | 47 | 20 | 1 | | 19 | 2 | 61 | 48 | | 34 | 51 | 16 | 29 | | 15 | 30 | 33 | 52 |
| 39 | 26 | 9 | 56 | | 10 | 55 | 40 | 25 | | 23 | 42 | 57 | 8 | | 58 | 7 | 24 | 41 |

186. Awani Kumar 2006    D1= 102    D2= 166    D3= 94    D4= 158

| 18 | 5 | 44 | 63 | | 45 | 58 | 23 | 4 | | 12 | 31 | 50 | 37 | | 55 | 36 | 13 | 26 |
|---|---|---|---|---|---|---|---|---|---|---|---|---|---|---|---|---|---|---|
| 11 | 32 | 49 | 38 | | 56 | 35 | 14 | 25 | | 17 | 6 | 43 | 64 | | 46 | 57 | 24 | 3 |
| 48 | 59 | 22 | 1 | | 19 | 8 | 41 | 62 | | 54 | 33 | 16 | 27 | | 9 | 30 | 51 | 40 |
| 53 | 34 | 15 | 28 | | 10 | 29 | 52 | 39 | | 47 | 60 | 21 | 2 | | 20 | 7 | 42 | 61 |

187. Guenter Stertenbrink 2003    D1= 130    D2= 130    D3= 130    D4= 130

| 18 | 5 | 44 | 63 | | 45 | 58 | 23 | 4 | | 12 | 31 | 50 | 37 | | 55 | 36 | 13 | 26 |
|---|---|---|---|---|---|---|---|---|---|---|---|---|---|---|---|---|---|---|
| 43 | 64 | 17 | 6 | | 24 | 3 | 46 | 57 | | 49 | 38 | 11 | 32 | | 14 | 25 | 56 | 35 |
| 16 | 59 | 22 | 33 | | 19 | 40 | 9 | 62 | | 54 | 1 | 48 | 27 | | 41 | 30 | 51 | 8 |
| 53 | 2 | 47 | 28 | | 42 | 29 | 52 | 7 | | 15 | 60 | 21 | 34 | | 20 | 39 | 10 | 61 |

188. Awani Kumar 2017    D1= 130    D2= 130    D3= 130    D4= 130

| 18 | 5 | 44 | 63 | | 45 | 58 | 23 | 4 | | 12 | 31 | 50 | 37 | | 55 | 36 | 13 | 26 |
|---|---|---|---|---|---|---|---|---|---|---|---|---|---|---|---|---|---|---|
| 43 | 64 | 17 | 6 | | 24 | 3 | 46 | 57 | | 49 | 38 | 11 | 32 | | 14 | 25 | 56 | 35 |
| 48 | 59 | 22 | 1 | | 19 | 8 | 41 | 62 | | 54 | 33 | 16 | 27 | | 9 | 30 | 51 | 40 |
| 21 | 2 | 47 | 60 | | 42 | 61 | 20 | 7 | | 15 | 28 | 53 | 34 | | 52 | 39 | 10 | 29 |

189. Francis Gaspalou 2009    D1= 66    D2= 194    D3= 66    D4= 194

| 18 | 7 | 48 | 57 | | 45 | 60 | 19 | 6 | | 40 | 49 | 26 | 15 | | 27 | 14 | 37 | 52 |
|---|---|---|---|---|---|---|---|---|---|---|---|---|---|---|---|---|---|---|
| 41 | 64 | 23 | 2 | | 22 | 3 | 44 | 61 | | 31 | 10 | 33 | 56 | | 36 | 53 | 30 | 11 |
| 8 | 17 | 58 | 47 | | 59 | 46 | 5 | 20 | | 50 | 39 | 16 | 25 | | 13 | 28 | 51 | 38 |
| 63 | 42 | 1 | 24 | | 4 | 21 | 62 | 43 | | 9 | 32 | 55 | 34 | | 54 | 35 | 12 | 29 |

190. Francis Gaspalou 2009    D1= 66    D2= 194    D3= 194    D4= 66

| 18 | 9 | 40 | 63 | | 45 | 54 | 27 | 4 | | 8 | 31 | 50 | 41 | | 59 | 36 | 13 | 22 |
|---|---|---|---|---|---|---|---|---|---|---|---|---|---|---|---|---|---|---|
| 39 | 64 | 17 | 10 | | 28 | 3 | 46 | 53 | | 49 | 42 | 7 | 32 | | 14 | 21 | 60 | 35 |
| 16 | 55 | 26 | 33 | | 19 | 44 | 5 | 62 | | 58 | 1 | 48 | 23 | | 37 | 30 | 51 | 12 |
| 57 | 2 | 47 | 24 | | 38 | 29 | 52 | 11 | | 15 | 56 | 25 | 34 | | 20 | 43 | 6 | 61 |

191. Awani Kumar 2017    D1= 130    D2= 130    D3= 130    D4= 130

| 18 | 13 | 40 | 59 | | 41 | 54 | 31 | 4 | | 8 | 27 | 50 | 45 | | 63 | 36 | 9 | 22 |
|---|---|---|---|---|---|---|---|---|---|---|---|---|---|---|---|---|---|---|
| 39 | 60 | 17 | 14 | | 32 | 3 | 42 | 53 | | 49 | 46 | 7 | 28 | | 10 | 21 | 64 | 35 |
| 48 | 19 | 58 | 5 | | 55 | 12 | 33 | 30 | | 26 | 37 | 16 | 51 | | 1 | 62 | 23 | 44 |
| 25 | 38 | 15 | 52 | | 2 | 61 | 24 | 43 | | 47 | 20 | 57 | 6 | | 56 | 11 | 34 | 29 |

192. Awani Kumar 2017    D1= 66    D2= 194    D3= 66    D4= 194

| 18 | 13 | 44 | 55 | | 45 | 50 | 23 | 12 | | 36 | 63 | 26 | 5 | | 31 | 4 | 37 | 58 |
|---|---|---|---|---|---|---|---|---|---|---|---|---|---|---|---|---|---|---|
| 35 | 64 | 17 | 14 | | 32 | 3 | 46 | 49 | | 25 | 6 | 43 | 56 | | 38 | 57 | 24 | 11 |
| 16 | 51 | 54 | 9 | | 19 | 48 | 41 | 22 | | 62 | 1 | 8 | 59 | | 33 | 30 | 27 | 40 |
| 61 | 2 | 15 | 52 | | 34 | 29 | 20 | 47 | | 7 | 60 | 53 | 10 | | 28 | 39 | 42 | 21 |

193. Awani Kumar 2017    D1= 50    D2= 130    D3= 210    D4= 130

| 18 | 13 | 44 | 55 | | 43 | 56 | 17 | 14 | | 16 | 19 | 54 | 41 | | 53 | 42 | 15 | 20 |
|---|---|---|---|---|---|---|---|---|---|---|---|---|---|---|---|---|---|---|
| 37 | 58 | 31 | 4 | | 32 | 3 | 38 | 57 | | 59 | 40 | 1 | 30 | | 2 | 29 | 60 | 39 |
| 12 | 23 | 50 | 45 | | 49 | 46 | 11 | 24 | | 22 | 9 | 48 | 51 | | 47 | 52 | 21 | 10 |
| 63 | 36 | 5 | 26 | | 6 | 25 | 64 | 35 | | 33 | 62 | 27 | 8 | | 28 | 7 | 34 | 61 |

194. Guenter Stertenbrink 2003    D1= 130    D2= 130    D3= 130    D4= 130

| 18 | 13 | 44 | 55 | | 43 | 56 | 17 | 14 | | 48 | 19 | 54 | 9 | | 21 | 42 | 15 | 52 |
|---|---|---|---|---|---|---|---|---|---|---|---|---|---|---|---|---|---|---|
| 37 | 58 | 31 | 4 | | 32 | 3 | 38 | 57 | | 59 | 8 | 33 | 30 | | 2 | 61 | 28 | 39 |
| 12 | 23 | 50 | 45 | | 49 | 46 | 11 | 24 | | 22 | 41 | 16 | 51 | | 47 | 20 | 53 | 10 |
| 63 | 36 | 5 | 26 | | 6 | 25 | 64 | 35 | | 1 | 62 | 27 | 40 | | 60 | 7 | 34 | 29 |

195. Awani Kumar 2017    D1= 66    D2= 194    D3= 194    D4= 66

| 18 | 13 | 44 | 55 | | 45 | 50 | 23 | 12 | | 36 | 63 | 26 | 5 | | 31 | 4 | 37 | 58 |
|---|---|---|---|---|---|---|---|---|---|---|---|---|---|---|---|---|---|---|
| 43 | 56 | 25 | 6 | | 24 | 11 | 38 | 57 | | 17 | 14 | 35 | 64 | | 46 | 49 | 32 | 3 |
| 16 | 51 | 54 | 9 | | 19 | 48 | 41 | 22 | | 62 | 1 | 8 | 59 | | 33 | 30 | 27 | 40 |
| 53 | 10 | 7 | 60 | | 42 | 21 | 28 | 39 | | 15 | 52 | 61 | 2 | | 20 | 47 | 34 | 29 |

196. Awani Kumar 2017    D1= 66    D2= 114    D3= 194    D4= 146

| 18 | 39 | 48 | 25 | | 45 | 28 | 19 | 38 | | 40 | 17 | 26 | 47 | | 27 | 46 | 37 | 20 |
|---|---|---|---|---|---|---|---|---|---|---|---|---|---|---|---|---|---|---|
| 41 | 32 | 23 | 34 | | 22 | 35 | 44 | 29 | | 31 | 42 | 33 | 24 | | 36 | 21 | 30 | 43 |
| 8 | 49 | 58 | 15 | | 59 | 14 | 5 | 52 | | 50 | 7 | 16 | 57 | | 13 | 60 | 51 | 6 |
| 63 | 10 | 1 | 56 | | 4 | 53 | 62 | 11 | | 9 | 64 | 55 | 2 | | 54 | 3 | 12 | 61 |

197. Francis Gaspalou 2009    D1= 130    D2= 130    D3= 130    D4= 130

| 20 | 1  | 62 | 47 |
|----|----|----|----|
| 9  | 28 | 39 | 54 |
| 64 | 45 | 18 | 3  |
| 37 | 56 | 11 | 26 |

| 61 | 48 | 19 | 2  |
|----|----|----|----|
| 40 | 53 | 10 | 27 |
| 17 | 4  | 63 | 46 |
| 12 | 25 | 38 | 55 |

| 16 | 29 | 34 | 51 |
|----|----|----|----|
| 21 | 8  | 59 | 42 |
| 36 | 49 | 14 | 31 |
| 57 | 44 | 23 | 6  |

| 33 | 52 | 15 | 30 |
|----|----|----|----|
| 60 | 41 | 22 | 7  |
| 13 | 32 | 35 | 50 |
| 24 | 5  | 58 | 43 |

198. Awani Kumar 2006 — D1= 130 D2= 130 D3= 130 D4= 130

| 20 | 5  | 46 | 59 |
|----|----|----|----|
| 43 | 62 | 21 | 4  |
| 6  | 19 | 60 | 45 |
| 61 | 44 | 3  | 22 |

| 47 | 58 | 17 | 8  |
|----|----|----|----|
| 24 | 1  | 42 | 63 |
| 57 | 48 | 7  | 18 |
| 2  | 23 | 64 | 41 |

| 38 | 51 | 28 | 13 |
|----|----|----|----|
| 29 | 12 | 35 | 54 |
| 52 | 37 | 14 | 27 |
| 11 | 30 | 53 | 36 |

| 25 | 16 | 39 | 50 |
|----|----|----|----|
| 34 | 55 | 32 | 9  |
| 15 | 26 | 49 | 40 |
| 56 | 33 | 10 | 31 |

199. Francis Gaspalou 2009 — D1= 66 D2= 194 D3= 194 D4= 66

| 20 | 7  | 42 | 61 |
|----|----|----|----|
| 41 | 62 | 19 | 8  |
| 46 | 25 | 56 | 3  |
| 23 | 36 | 13 | 58 |

| 47 | 60 | 21 | 2  |
|----|----|----|----|
| 22 | 1  | 48 | 59 |
| 49 | 6  | 43 | 32 |
| 12 | 63 | 18 | 37 |

| 10 | 29 | 52 | 39 |
|----|----|----|----|
| 51 | 40 | 9  | 30 |
| 24 | 35 | 14 | 57 |
| 45 | 26 | 55 | 4  |

| 53 | 34 | 15 | 28 |
|----|----|----|----|
| 16 | 27 | 54 | 33 |
| 11 | 64 | 17 | 38 |
| 50 | 5  | 44 | 31 |

200. Awani Kumar 2017 — D1= 66 D2= 194 D3= 66 D4= 194

| 20 | 7  | 42 | 61 |
|----|----|----|----|
| 41 | 62 | 19 | 8  |
| 46 | 57 | 24 | 3  |
| 23 | 4  | 45 | 58 |

| 47 | 60 | 21 | 2  |
|----|----|----|----|
| 22 | 1  | 48 | 59 |
| 17 | 6  | 43 | 64 |
| 44 | 63 | 18 | 5  |

| 10 | 29 | 52 | 39 |
|----|----|----|----|
| 51 | 40 | 9  | 30 |
| 56 | 35 | 14 | 25 |
| 13 | 26 | 55 | 36 |

| 53 | 34 | 15 | 28 |
|----|----|----|----|
| 16 | 27 | 54 | 33 |
| 11 | 32 | 49 | 38 |
| 50 | 37 | 12 | 31 |

201. Francis Gaspalou 2009 — D1= 66 D2= 194 D3= 66 D4= 194

| 20 | 11 | 38 | 61 |
|----|----|----|----|
| 37 | 62 | 19 | 12 |
| 46 | 21 | 60 | 3  |
| 27 | 36 | 13 | 54 |

| 47 | 56 | 25 | 2  |
|----|----|----|----|
| 26 | 1  | 48 | 55 |
| 49 | 10 | 39 | 32 |
| 8  | 63 | 18 | 41 |

| 6  | 29 | 52 | 43 |
|----|----|----|----|
| 51 | 44 | 5  | 30 |
| 28 | 35 | 14 | 53 |
| 45 | 22 | 59 | 4  |

| 57 | 34 | 15 | 24 |
|----|----|----|----|
| 16 | 23 | 58 | 33 |
| 7  | 64 | 17 | 42 |
| 50 | 9  | 40 | 31 |

202. Awani Kumar 2017 — D1= 66 D2= 194 D3= 66 D4= 194

| 20 | 15 | 42 | 53 |
|----|----|----|----|
| 41 | 54 | 27 | 8  |
| 46 | 17 | 24 | 43 |
| 23 | 44 | 37 | 26 |

| 47 | 52 | 21 | 10 |
|----|----|----|----|
| 22 | 9  | 40 | 59 |
| 49 | 14 | 11 | 56 |
| 12 | 55 | 58 | 5  |

| 34 | 61 | 28 | 7  |
|----|----|----|----|
| 19 | 16 | 33 | 62 |
| 32 | 35 | 38 | 25 |
| 45 | 18 | 31 | 36 |

| 29 | 2  | 39 | 60 |
|----|----|----|----|
| 48 | 51 | 30 | 1  |
| 3  | 64 | 57 | 6  |
| 50 | 13 | 4  | 63 |

203. Awani Kumar 2017 — D1= 130 D2= 178 D3= 130 D4= 82

| 21 | 2  | 47 | 60 |
|----|----|----|----|
| 16 | 27 | 54 | 33 |
| 43 | 64 | 17 | 6  |
| 50 | 37 | 12 | 31 |

| 42 | 61 | 20 | 7  |
|----|----|----|----|
| 51 | 40 | 9  | 30 |
| 24 | 3  | 46 | 57 |
| 13 | 26 | 55 | 36 |

| 15 | 28 | 53 | 34 |
|----|----|----|----|
| 22 | 1  | 48 | 59 |
| 49 | 38 | 11 | 32 |
| 44 | 63 | 18 | 5  |

| 52 | 39 | 10 | 29 |
|----|----|----|----|
| 41 | 62 | 19 | 8  |
| 14 | 25 | 56 | 35 |
| 23 | 4  | 45 | 58 |

204. Guenter Stertenbrink 2003 — D1= 130 D2= 130 D3= 130 D4= 130

| 21 | 10 | 47 | 52 |
|----|----|----|----|
| 34 | 61 | 28 | 7  |
| 15 | 20 | 53 | 42 |
| 60 | 39 | 2  | 29 |

| 48 | 51 | 22 | 9  |
|----|----|----|----|
| 27 | 8  | 33 | 62 |
| 54 | 41 | 16 | 19 |
| 1  | 30 | 59 | 40 |

| 11 | 24 | 49 | 46 |
|----|----|----|----|
| 64 | 35 | 6  | 25 |
| 17 | 14 | 43 | 56 |
| 38 | 57 | 32 | 3  |

| 50 | 45 | 12 | 23 |
|----|----|----|----|
| 5  | 26 | 63 | 36 |
| 44 | 55 | 18 | 13 |
| 31 | 4  | 37 | 58 |

205. Guenter Stertenbrink 2003 — D1= 130 D2= 130 D3= 130 D4= 130

| 21 | 10 | 47 | 52 |
|----|----|----|----|
| 34 | 61 | 28 | 7  |
| 15 | 20 | 53 | 42 |
| 60 | 39 | 2  | 29 |

| 48 | 51 | 22 | 9  |
|----|----|----|----|
| 27 | 8  | 33 | 62 |
| 54 | 41 | 16 | 19 |
| 1  | 30 | 59 | 40 |

| 11 | 56 | 17 | 46 |
|----|----|----|----|
| 32 | 35 | 6  | 57 |
| 49 | 14 | 43 | 24 |
| 38 | 25 | 64 | 3  |

| 50 | 13 | 44 | 23 |
|----|----|----|----|
| 37 | 26 | 63 | 4  |
| 12 | 55 | 18 | 45 |
| 31 | 36 | 5  | 58 |

206. Awani Kumar 2017 — D1= 130 D2= 130 D3= 130 D4= 130

| | | | | | | | | | | | | | | | |
|---|---|---|---|---|---|---|---|---|---|---|---|---|---|---|---|
| 21 | 10 | 47 | 52 | 48 | 51 | 30 | 1 | 11 | 56 | 49 | 14 | 50 | 13 | 4 | 63 |
| 42 | 53 | 20 | 15 | 19 | 16 | 33 | 62 | 24 | 43 | 46 | 17 | 45 | 18 | 31 | 36 |
| 39 | 60 | 29 | 2 | 22 | 9 | 40 | 59 | 57 | 6 | 3 | 64 | 12 | 55 | 58 | 5 |
| 28 | 7 | 34 | 61 | 41 | 54 | 27 | 8 | 38 | 25 | 32 | 35 | 23 | 44 | 37 | 26 |

207. Awani Kumar 2017    D1= 66   D2= 114   D3= 146   D4= 194

| | | | | | | | | | | | | | | | |
|---|---|---|---|---|---|---|---|---|---|---|---|---|---|---|---|
| 22 | 35 | 44 | 29 | 45 | 28 | 19 | 38 | 4 | 53 | 62 | 11 | 59 | 14 | 5 | 52 |
| 41 | 32 | 23 | 34 | 18 | 39 | 48 | 25 | 63 | 10 | 1 | 56 | 8 | 49 | 58 | 15 |
| 36 | 21 | 30 | 43 | 27 | 46 | 37 | 20 | 54 | 3 | 12 | 61 | 13 | 60 | 51 | 6 |
| 31 | 42 | 33 | 24 | 40 | 17 | 26 | 47 | 9 | 64 | 55 | 2 | 50 | 7 | 16 | 57 |

208. Francis Gaspalou 2009    D1= 130   D2= 130   D3= 130   D4= 130

| | | | | | | | | | | | | | | | |
|---|---|---|---|---|---|---|---|---|---|---|---|---|---|---|---|
| 23 | 2 | 59 | 46 | 58 | 47 | 22 | 3 | 7 | 62 | 43 | 18 | 42 | 19 | 6 | 63 |
| 10 | 51 | 38 | 31 | 39 | 30 | 11 | 50 | 26 | 15 | 54 | 35 | 55 | 34 | 27 | 14 |
| 57 | 48 | 21 | 4 | 24 | 1 | 60 | 45 | 41 | 20 | 5 | 64 | 8 | 61 | 44 | 17 |
| 40 | 29 | 12 | 49 | 9 | 52 | 37 | 32 | 56 | 33 | 28 | 13 | 25 | 16 | 53 | 36 |

209. Awani Kumar 2006    D1= 94   D2= 102   D3= 158   D4= 166

| | | | | | | | | | | | | | | | |
|---|---|---|---|---|---|---|---|---|---|---|---|---|---|---|---|
| 25 | 6 | 47 | 52 | 48 | 51 | 26 | 5 | 7 | 60 | 17 | 46 | 50 | 13 | 40 | 27 |
| 34 | 61 | 24 | 11 | 23 | 12 | 33 | 62 | 32 | 35 | 10 | 53 | 41 | 22 | 63 | 4 |
| 15 | 20 | 57 | 38 | 58 | 37 | 16 | 19 | 49 | 14 | 39 | 28 | 8 | 59 | 18 | 45 |
| 56 | 43 | 2 | 29 | 1 | 30 | 55 | 44 | 42 | 21 | 64 | 3 | 31 | 36 | 9 | 54 |

210. Awani Kumar 2017    D1= 130   D2= 130   D3= 130   D4= 130

| | | | | | | | | | | | | | | | |
|---|---|---|---|---|---|---|---|---|---|---|---|---|---|---|---|
| 34 | 7 | 32 | 57 | 29 | 60 | 39 | 2 | 24 | 33 | 42 | 31 | 43 | 30 | 17 | 40 |
| 25 | 64 | 35 | 6 | 38 | 3 | 28 | 61 | 47 | 26 | 21 | 36 | 20 | 37 | 46 | 27 |
| 8 | 49 | 58 | 15 | 59 | 14 | 1 | 56 | 50 | 23 | 16 | 41 | 13 | 44 | 55 | 18 |
| 63 | 10 | 5 | 52 | 4 | 53 | 62 | 11 | 9 | 48 | 51 | 22 | 54 | 19 | 12 | 45 |

211. Awani Kumar 2012    D1= 98   D2= 162   D3= 138   D4= 122

| | | | | | | | | | | | | | | | |
|---|---|---|---|---|---|---|---|---|---|---|---|---|---|---|---|
| 34 | 15 | 30 | 51 | 27 | 54 | 39 | 10 | 14 | 35 | 50 | 31 | 55 | 26 | 11 | 38 |
| 19 | 62 | 47 | 2 | 42 | 7 | 22 | 59 | 63 | 18 | 3 | 46 | 6 | 43 | 58 | 23 |
| 16 | 33 | 52 | 29 | 53 | 28 | 9 | 40 | 36 | 13 | 32 | 49 | 25 | 56 | 37 | 12 |
| 61 | 20 | 1 | 48 | 8 | 41 | 60 | 21 | 17 | 64 | 45 | 4 | 44 | 5 | 24 | 57 |

212. Awani Kumar 2006    D1= 130   D2= 130   D3= 130   D4= 130

| | | | | | | | | | | | | | | | |
|---|---|---|---|---|---|---|---|---|---|---|---|---|---|---|---|
| 34 | 23 | 32 | 41 | 29 | 44 | 39 | 18 | 24 | 49 | 42 | 15 | 43 | 14 | 17 | 56 |
| 25 | 48 | 35 | 22 | 38 | 19 | 28 | 45 | 47 | 10 | 21 | 52 | 20 | 53 | 46 | 11 |
| 8 | 33 | 58 | 31 | 59 | 30 | 1 | 40 | 50 | 7 | 16 | 57 | 13 | 60 | 55 | 2 |
| 63 | 26 | 5 | 36 | 4 | 37 | 62 | 27 | 9 | 64 | 51 | 6 | 54 | 3 | 12 | 61 |

213. Awani Kumar 2012    D1= 130   D2= 130   D3= 170   D4= 90

| | | | | | | | | | | | | | | | |
|---|---|---|---|---|---|---|---|---|---|---|---|---|---|---|---|
| 36 | 13 | 32 | 49 | 25 | 56 | 37 | 12 | 16 | 33 | 52 | 29 | 53 | 28 | 9 | 40 |
| 17 | 64 | 45 | 4 | 44 | 5 | 24 | 57 | 61 | 20 | 1 | 48 | 8 | 41 | 60 | 21 |
| 14 | 35 | 50 | 31 | 55 | 26 | 11 | 38 | 34 | 15 | 30 | 51 | 27 | 54 | 39 | 10 |
| 63 | 18 | 3 | 46 | 6 | 43 | 58 | 23 | 19 | 62 | 47 | 2 | 42 | 7 | 22 | 59 |

214. Awani Kumar 2006    D1= 130   D2= 130   D3= 130   D4= 130

| | | | | | | | | | | | | | | | |
|---|---|---|---|---|---|---|---|---|---|---|---|---|---|---|---|
| 36 | 21 | 26 | 47 | 31 | 42 | 33 | 24 | 22 | 51 | 48 | 9 | 41 | 16 | 23 | 50 |
| 27 | 46 | 37 | 20 | 40 | 17 | 30 | 43 | 45 | 12 | 19 | 54 | 18 | 55 | 44 | 13 |
| 6 | 35 | 64 | 25 | 57 | 32 | 7 | 34 | 52 | 5 | 10 | 63 | 15 | 58 | 49 | 8 |
| 61 | 28 | 3 | 38 | 2 | 39 | 60 | 29 | 11 | 62 | 53 | 4 | 56 | 1 | 14 | 59 |

215. Awani Kumar 2012    D1= 122   D2= 138   D3= 162   D4= 98

| | | | | | | | | | | | | | | | |
|---|---|---|---|---|---|---|---|---|---|---|---|---|---|---|---|
| 38 | 19 | 28 | 45 | 29 | 44 | 39 | 18 | 4 | 37 | 62 | 27 | 59 | 30 | 1 | 40 |
| 25 | 48 | 35 | 22 | 34 | 23 | 32 | 41 | 63 | 26 | 5 | 36 | 8 | 33 | 58 | 31 |
| 20 | 53 | 46 | 11 | 43 | 14 | 17 | 56 | 54 | 3 | 12 | 61 | 13 | 60 | 55 | 2 |
| 47 | 10 | 21 | 52 | 24 | 49 | 42 | 15 | 9 | 64 | 51 | 6 | 50 | 7 | 16 | 57 |

216. Awani Kumar 2012    D1= 130   D2= 130   D3= 106   D4= 154